\documentclass{article}
\usepackage[utf8]{inputenc}

\usepackage{amsmath, amssymb, amsthm, stackrel,mathtools}
\usepackage{bbm}
\usepackage{authblk}
\usepackage{titlesec, tikz,subcaption,pgf}
\usepackage{tikz-cd}
\usepackage{hyperref}
\usepackage{xypic}
\usetikzlibrary{arrows}
\usetikzlibrary{knots}
\usetikzlibrary{decorations.markings}

\newcommand{\rar}{\rightarrow}
\newcommand{\tr}{\textnormal{tr}}

\newcommand{\cat}[1]{\mathcal{#1}}
\newcommand{\Z}{\mathbb{Z}}
\newcommand{\Ad}{\textnormal{Ad}}
\newcommand{\Aut}{\textnormal{Aut}}
\newcommand{\id}{1}

\newcommand{\End}{\textnormal{End}}
\newcommand{\Id}{\textnormal{Id}}
\newcommand{\ev}{\textnormal{ev}}
\newcommand{\coev}{\textnormal{coev}}
\newcommand{\Vect}{\mathbf{Vect}}
\newcommand{\qdim}{\textnormal{qdim}}

\newcommand{\tn}[1]{\textnormal{#1}}

\newcommand{\ns}[1]{#1_\textrm{NS}}

\newcommand{\cns}{\ns{\cat{C}}(d)}

\newcommand{\Rep}{\mathbf{Rep}}
\newcommand{\Fun}{\textnormal{Fun}}

\newcommand{\dcentcat}[1]{\cat{Z}(\cat{#1})}
\newcommand{\repsu}{\Rep(\widehat{\mathfrak{su}(2)}_{d-2})}
\newcommand{\pic}{\tn{Pic}}
\newcommand{\C}{\mathbb{C}}

\newcommand{\mi}{\tn{-}}
\newcommand{\vd}{\cat{V}_d}
\newcommand{\pd}{\mathcal{P}_d}
\newcommand{\Zd}{\underline{\Z_d}}

\newtheorem{thma}{Theorem}

\newtheorem{thm}{Theorem}

\newtheorem{cor}[thm]{Corollary}
\newtheorem{lem}[thm]{Lemma}
\newtheorem{prop}[thm]{Proposition}
\newtheorem{facts}[thm]{Facts}
\theoremstyle{definition}

\newtheorem{df}[thm]{Definition}
\newtheorem{rmk}[thm]{Remark}

\newtheorem{notation}[thm]{Notation}

\title{Landau-Ginzburg/Conformal Field Theory Correspondence for $x^d$ and Module Tensor Categories}
\author[1]{Ana Ros Camacho}
\author[2]{Thomas A. Wasserman}
\affil[1]{School of Mathematics, Cardiff University}
\affil[2]{Lincoln College, University of Oxford}
\date{\today}

\begin{document}

\maketitle
\begin{abstract}
    The Landau-Ginzburg/Conformal Field Theory correspondence predicts tensor equivalences between categories of matrix factorisations of certain polynomials and categories associated to the $N=2$ supersymmetric conformal field theories. We realise this correspondence for $x^d$ for any $d$, where previous results were limited to odd $d$. 
    
    Our proof uses the fact that both sides of the correspondence carry the structure of module tensor categories over the category of $\Z_d$-graded vector spaces equipped with a non-degenerate braiding. This allows us to describe the CFT side as generated by a single object, and use this to efficiently provide a functor realising the tensor equivalence.
\end{abstract}

\section{Introduction}

We establish the Landau-Ginzburg/Conformal Field Theory (LG/CFT) correspondence for a family of self-defects of the potential $x^d$ for any $d>2$. This generalises the main theorem of \cite{Davydov2018} to even $d$. The LG/CFT correspondence appeared in the physics literature in the late 80s and early 90s, see e.g. \cite{Brunner2006,Howe1989a,Kastor1989,Fermi1989,Vafa1989a}. It comes from the observation that every $N=2$-supersymmetric Landau-Ginzburg model characterized by a potential $W$ at the infrared renormalisation group fixed point is a conformal field theory of central charge $c_W$. At the level of boundary conditions, this correspondence predicts equivalences of $\C$-linear categories between (on the CFT-side) categories of representations of vertex operator algebras (VOAs) and (on the LG-side) categories of matrix factorisations of $W$. At the level of self-defects, these equivalences should be tensor. There have been some recent efforts towards providing a mathematical framework for this correspondence, see e.g. \cite{Carqueville2016b,Davydov2018,RosCamacho2020}.

The only instance of LG/CFT realised as a tensor equivalence is that from \cite{Davydov2018}. Here, the first author and collaborators provided a tensor equivalence between the Neveu-Schwarz part $\mathcal{C}_{\mathrm{NS}} \left( d \right)$ of the category $\cat{C}(d)$ of representations of the super-VOA of the unitary $N=2$ minimal model with central charge $c=3 \frac{(d-2)}{d}$ and the permutation type subcategory $\mathcal{P}_d$ of the category of matrix factorisations of $x^d-y^d$, for odd $d$. However, independent of $d$'s parity the fusion rules for $\mathcal{P}_d$ (first described in \cite{Brunner2007}, fully established in \cite{Davydov2018}) agree with those of $\cns$. We prove that this agreement can be promoted to a tensor functor for any $d$:

\begin{thma}\label{maintheoremboring}
Let $d>2$. There is a tensor equivalence $\cns \cong \cat{P}_d$.
\end{thma}

Our proof strategy provides a fresh perspective on the structure of the categories involved. We regard them as $\vd$-module tensor categories (cf. \cite{Henriques2016}), where $\vd$ is the category of $\Z_d$-graded vector spaces, equipped with a non-degenerate braiding and trivial associator (Definition \ref{vddef}). These are (in particular) tensor categories $\cat{A}$ equipped with a central functor $\vd \rar \dcentcat{A}$. As $\vd \subset\cns$ as a braided subcategory, $\cns$ is naturally a $\vd$-module tensor category. This allows us to in effect quotient by $\vd$. This looks like:

\begin{thma}
The category $\cns$ is equivalent as $\cat{V}_d$-module tensor category to the pivotal free $\cat{V}_d$-module fusion category generated by a single object that is self-dual of charge $\mi 1$ and quantum dimension $2 \cos(\frac{\pi}{d})$.
\end{thma}
This is Theorem \ref{cnsgenerated} in the main text. Intuitively, this means that $\cns$ is generated by a single object that is self dual up to the action of $\mi 1 \in \vd$, and by freely acting by $\vd$. This is an instance of a generalisation of the well-known construction of fusion categories generated by a single self-dual object as a quotient of Temperley-Lieb categories. In this paper, we spell out the generalisation we need in the specific case of $\cns$. Functors out of a category generated by a single self-dual object are easy to describe: one just needs to identify a viable target for the object and its duality data.

On the LG-side, the $\cat{V}_d$-module tensor category structure comes from the action of automorphisms of $\C[x]$ fixing $x^d$. Giving a trivialisation of the action of $\Z_d$ on $\mathcal{P}_d$ that is consistent across $\cat{P}_d$ and compatible with the tensor structure yields (Theorem \ref{pdisvdmod} in the main text):

\begin{thma}
    The category $\cat{P}_d$ can be given the structure of a spherical $\vd$-module tensor category. 
\end{thma}

To establish tensor equivalence between $\cns$ and $\cat{P}_d$ we provide a self-dual object of charge $\mi 1$ in $\cat{P}_d$, giving rise to a tensor functor $\cns\rar \pd$. Standard arguments show that this is an equivalence:

\begin{thma}
The categories $\cns$ and $\cat{P}_d$ are equivalent as spherical free $\vd$-module tensor categories.
\end{thma}

This appears as Theorem \ref{badassmaintheorem}. Theorem \ref{maintheoremboring} is a direct corollary of this.

This raises interesting questions. Regarding our generalisation of the Temperley-Lieb construction, in \cite{Henriques2016a} another, more powerful, generalisation of this construction is given. However, our case falls slightly outside the applicability of their results. While it is beyond the scope of this work, we expect that our case can be fitted into their framework. 

On the topic of the LG/CFT correspondence, we observe the following. In the physics context the $\Z_d$-action arises as a symmetry of the Lagrangian defining the LG-model, and the associated chiral ring. This symmetry is preserved under the renormalization group flow, so the corresponding CFT would also have this symmetry group. Mathematically, this translates to an action by the automorphisms $G$ of the potential $W$ on the category of matrix (bi-)factorisations for $W$. We expect that this gives rise to the structure of a $\Vect[G]$-module category on the category of matrix factorisation for $W$, that in turn can be lifted to a module tensor category structure on a subcategory, after picking an appropriate quadratic form on $G$. Identifying the associated $\Vect[G]$-module category structure should give another hand-hold to establishing the LG/CFT correspondence in other examples, we will explore this in future work.

\subsubsection*{Outline}
The paper is organized as follows. In Section \ref{Section2} we review the category $\cns$ and some basics on Temperley-Lieb categories. We then describe the structure of $\cns$ as a $\mathcal{V}_d$-module tensor category. In Section \ref{Section3} we provide a brief overview of matrix factorisations and introduce the category $\pd$ of permutation type matrix factorisations. We then equip this category with the structure of a $\vd$-module tensor category. Finally, in Section \ref{Section4} we state and prove our main result.

\subsection{Preliminaries}
\subsubsection{Notation and conventions}
\label{notationsection}
We will work over $\mathbb{C}$. We fix the primitive $d$-th root of unity $\eta=e^{\frac{2\pi i}{d}}$ throughout. 

When tensoring two objects in a category we will usually omit (especially in diagrams) the $\otimes$ symbol and for tensor products of morphisms, we will use the notation $\otimes=\cdot$. We will write $\id_c$ for the identity on an object $c$, and will often omit identity morphisms tensored with any other morphisms if clear from the context.

We will use string diagram calculus for braided fusion categories to streamline our arguments. Readers unfamiliar with this are referred to e.g. \cite{Joyal1993,Bartlett2015}. As is customary, unitors and associators are suppressed in the string diagrams.

\subsubsection{The category $\vd$}
We remind the reader of the classical result \cite{Joyal1986} that braidings on monoidal categories with a finite group as their set of objects are completely determined by a quadratic form on the group. This allows us to define:

\begin{df}\label{vddef}
The pointed ribbon fusion category $\vd$ is the spherical fusion category of $\Z_d$-graded vector spaces equipped with the braiding and associators induced by the quadratic form $q_d(k)=\eta^{\frac{2\pi i}{d}k^2}$.
\end{df}

\begin{lem}\label{jscomputation}
The fusion category $\cat{V}_d$ has trivial associators. 
\end{lem}
\begin{proof}
Following the proof of \cite[Theorem 12]{Joyal1986}, we compute that the associators for the objects $l,m,n\in \Z_d$ on $\vd$ determined by the quadratic form are given by $\id_{l+m+n}$ times
$$
h(l,m,n)= \begin{cases}1 &\mbox{ for }m+n< d \\ (q_d(1)^{d})^l=1^l=1 & \mbox{else.} \end{cases}
$$
These are indeed trivial.
\end{proof}

\section{The conformal field theoretical side}\label{Section2}
In this section we set up the CFT side of the correspondence. We will briefly recall the categorical structure of the $N=2$ unitary minimal model with central charge $c=3\frac{(d-2)}{d}$. After this we discuss Temperley-Lieb categories, which will serve as a template for our description of $\cns$ as generated by a single object. This first part of this section is essentially a summary of \cite[Sections 2.1 and 2.3]{Davydov2018}. We then move on to our result describing $\cns$ as a $\vd$-module tensor category, first recalling some basics of module tensor categories, and then proving that $\cns$ can be viewed as a module tensor category generated by a single object.

\subsection{The category $\mathcal{C}_{\mathrm{NS}} \left( d \right)$}

\subsubsection{Coset construction of $\cns$}
Let $V$ be the vertex operator superalgebra corresponding to the $N=2$ unitary minimal model with central charge $c=3\frac{(d-2)}{d}$ with $d \in \mathbb{Z}_{\geq 2}$ \cite{Adamovic1999}. The degree-zero part $V_0$ of this superVOA can be identified with the coset $\frac{\widehat{\mathfrak{su}\left( 2 \right)}_{d-2} \oplus \widehat{\mathfrak{u \left( 1 \right)}}_4}{\widehat{\mathfrak{u \left( 1 \right)}}_{2d}}$ \cite{Carpi2015,DiVecchia1986}. The category of representations of this part of the VOA can be described as follows \cite{Frohlich2004}. Consider the category
\begin{equation*}
\mathcal{C}\left( d \right):= \mathbf{Rep}\left( \widehat{\mathfrak{su} \left(2 \right)}_{d-2} \right) \boxtimes \mathbf{Rep}\left( \widehat{\mathfrak{u} \left(1 \right)}_{2d} \right)^{\mathrm{op}} \boxtimes \mathbf{Rep}\left( \widehat{\mathfrak{u} \left(1 \right)}_{4} \right).
\end{equation*}
Here $\repsu$ is the category of integrable highest weight representations of the affine Lie algebra $\mathfrak{su} \left( 2 \right)$ at the level $d-2$, and $\mathbf{Rep}\left( \widehat{\mathfrak{u} \left(1 \right)}_{2d} \right)$ is the category of representations of the vertex operator algebra for $\mathfrak{u}\left( 1 \right)$ rationally extended by two fields of weight $d$. This is the pointed fusion category with simple objects represented by $\Z_{2d}$ and braiding and associators induced from the quadratic form on $\Z_{2d}$ given by
    \begin{equation}\label{quadformzd}
    	\widetilde{q}_{2d}\colon r \mapsto e^{\frac{i\pi r^2}{2d}}.
    \end{equation}
The superscript $^{op}$ indicates taking the opposite braiding and ribbon twist.

Simples in $\mathcal{C} \left( d \right)$ are labeled by triples $\left[ l,r,s \right]$ with $l \in \lbrace 0,\ldots,d-2 \rbrace$, $m \in \mathbb{Z}_{2d}$ and $s \in \mathbb{Z}_4$. The category of representations of $V_0$ can be realized as the category of local modules over the object $A:=\left[ 0,0,0 \right] \oplus \left[ d-2,d,2 \right]$ in $\mathcal{C} \left( d \right)$. Because $\left[ d-2,d,2 \right]$ is invertible and acts freely, simple $A$-modules are free. That is, simple modules are of the form $A\otimes \left[ l,r,s \right]$, they are local when $l+m+s$ is even. 

We are interested in the Neveu-Schwarz subcategory $\cns$, this consists of those objects with $s$ even.
In particular, we can write any simple of $\cns$ as $\left[ l,r \right]:=A \otimes \left[ l,r,0 \right]$.

\subsubsection{Categorical structure}
We will work at the level of the braided fusion categorical structure of $\cns$. The fusion rules are
\begin{equation}\label{fusionincns}
    [l,r]\otimes [l',r']= \bigoplus\limits_{\nu=\vert l-l' \vert \mbox{ step } 2}^{\mathrm{min\left( l+l',2d-4-l-l' \right)}} [\nu, r+r'].
\end{equation}
The associators on the $\repsu$-factor are rather complicated, see \cite[Appendix A.2]{Carqueville2010}, and we will not work with them directly. $\cns$ is ribbon, hence spherical, and the quantum dimensions are:
$$
\qdim([l,r])=\frac{\zeta^{l+1}-\zeta^{\mi l \mi 1}}{\zeta-\zeta^{-1}},
$$
where $\zeta=e^{\frac{\pi i}{d}}$ (so that $\zeta^2=\eta$).
From the fusion rules we see that there is a fusion subcategory spanned by $[0,2k]$ for $k\in \Z_d$. This is a copy of $\vd$, its quadratic form is that from Equation \eqref{quadformzd} pulled back to $\Z_d$:
\begin{equation}\label{vtwodtovd}
\widetilde{q}_{2d}(2k)=e^{\frac{i\pi (2k)^2}{2d}}=q_d(k).    
\end{equation}

\subsubsection{Temperley-Lieb categories} \label{temperleysection}
We want to describe $\cns$ as category generated by a single object in the appropriate sense. To do this, we use a generalisation of the construction the classical construction of fusion categories generated by a single self-dual object using Temperley-Lieb categories. In this section we recall this classical construction, closely following the treatment in \cite{Yamagami2004}. To start, fix a complex number $\kappa$.

Recall that a \textit{Kauffman diagram of type $\left( m,n \right)$} is a planar curve connecting a line of $m$ points to a line of $n$ points, confined to the rectangle between these lines. Note that there are no such curves when $m$ and $n$ have different parity. We denote the set of isotopy classes Kauffman diagrams by $K_{m,n}$. Concatenating diagrams induces maps $K_{m,n}\times K_{n,r}\rar K_{m,r}$, juxtaposition gives maps $K_{m,n}\times K_{m',n'}\rar K_{m+m',n+n'}$.

\begin{df}
The \emph{Temperley-Lieb category $\mathcal{TL} \left({\kappa} \right)$} is the idempotent completion of the tensor category $T_\kappa$ with objects spanned by the non-negative integers. For objects $m$ and $n$, we set $T_\kappa\left( m,n \right)=\mathbb{C} \left[ K_{m,n} \right]$, composition of morphisms is induced from concatenation of diagrams. We impose the relation that a circle is equivalent to the empty diagram times $\kappa$. The tensor product is addition on the objects and induced from juxtapostion on the morphisms.
\end{df}

 This category is generated by the object corresponding to $1\in \mathbb{N}$. This object is self-dual, with evaluation and co-evaluation given by the cap- and cup-shaped Kauffman diagrams. Its quantum dimension is $\kappa$. This category is spherical by construction.
 
To find the idempotent completion, we need to examine the endomorphism algebras in $T_\kappa$. These algebras are the \textit{Temperley-Lieb algebras} $\tn{TL}_n \left( {\kappa} \right)$. Adding a strand on the right of Kauffman diagrams in $K_{n,n}$ gives a map $K_{n,n}\rar K_{n+1,n+1}$, which induces a map on the Temperley-Lieb algebras $\tn{TL}_n \left( \kappa \right) \hookrightarrow \tn{TL}_{n+1} \left( \kappa \right)$. This allows one to inductively construct a sequence of idempotents called the \textit{Jones-Wenzl projectors} $p_{n}\in \tn{TL}_{n} \left( \kappa \right)$. Pick $\zeta$ such that $\kappa=\zeta+ \zeta^{-1}$. The induction terminates when the so-called \textit{quantum number} 
$$
[n+1]_{\zeta}= \frac{\zeta^{n+1}+\zeta^{-n-1}}{\zeta-\zeta^{-1}}
$$
vanishes. When $\zeta=e^{\frac{i\pi}{d}}$, this happens when $n=d-1$, and we have that $\tr(p_{d-1})=[d]_{\zeta}=0$ where $\tr$ denotes the Markov trace. For further properties of the Wenzl-Jones idempotents, as well as of Temperley-Lieb categories, we refer to e.g. \cite{Yamagami2004,Goodman2003}.

These idempotents give rise to a sequence of subobjects $\langle\!\langle n \rangle\!\rangle\subset n$, and for each $n< d-1$ the rest of the direct summands of $n$ can be identified with $\langle\!\langle i \rangle\!\rangle\subset$ for $i<n$ of the same parity as $n$. This gives the fusion rules for this category, the fusion morphisms correspond to elements of $K_{n,i}$ sandwiched between $p_n$ and $p_i$.

The final well-defined projector $p_{d-1}$ is what is called a negligible morphism (see e.g. \cite[Section 1.3]{Goodman2003}). It in fact generates the tensor ideal of negligible morphisms \cite[Proposition 2.1]{Goodman2003}. The quotient by this ideal is a spherical fusion category. Mapping the object associated to $1$ to the object $[1]\in \repsu$ induces a tensor functor $T_\kappa\rar \repsu$, which factors through this quotient as the image of $p_{d-1}$ in $\repsu$ is zero. Using semi-simplicity one then sees that this functor is an equivalence identifying this quotient with the spherical fusion category underlying $\repsu$.

\subsection{$\mathcal{C}_{\tn{NS}}$ as a $\cat{V}_d$-module tensor category}
For $d$ odd the category $\mathcal{C}_{\mathrm{NS}} \left( d \right)$ splits and is equivalent as a braided fusion category to the Deligne tensor product $\frac{\mathcal{TL}_{2 \mathrm{cos} (\frac{\pi}{d})}}{\langle p_{d-1} \rangle} \boxtimes \mathcal{V}_d$, by \cite[Proposition 2.1]{Davydov2018}.

For even $d$, this does not hold. What fails is that one is not able to find a self-dual object to generate the Temperley-Lieb factor. In this section we will show this can be circumvented by taking into account that $\vd$ acts on $\cns$.

\subsubsection{Module tensor categories}
The following definition (see \cite[Definition 3.1]{Henriques2016a} and \cite[Section 3.2]{Henriques2016} for more details) allows us to capture the action of $\vd$ on $\cns$:

\begin{df}
Let $\mathcal{D}$ be a braided tensor category. A \textit{module tensor category} over $\mathcal{D}$ is a tensor category $\mathcal{M}$ together with a braided tensor functor $\Phi^\mathcal{Z} \colon \mathcal{D} \to \mathcal{Z} \left( \mathcal{M} \right)$ from $\mathcal{D}$ to the Drinfeld center of $\mathcal{M}$. We will call $\cat{M}$ a module \emph{fusion} category if both $\cat{D}$ and $\mathcal{M}$ are fusion, and similarly \emph{pivotal} (or \emph{spherical}) if $\cat{D}$, $\mathcal{M}$ and $\Phi^\mathcal{Z}$ are.

Write $\Phi$ for the composite of $\Phi^\cat{Z}$ with the forgetful functor $\dcentcat{M}\to \mathcal{M}$. Then $\mathcal{M}$ is called \emph{free} if for any simple $m \in \cat{M}$ and simples $v,v'\in \cat{D}$ we have $\Phi(v) m \cong \Phi(v') m$ if and only if $v\cong v'$.
\end{df}

Equation \eqref{vtwodtovd} implies that $\cns$ is a $\cat{V}_d$-module tensor category: it is a tensor category equipped with a braided functor $\Phi: \cat{V}_d \hookrightarrow \cat{Z}(\cns)$, which in this case factors through $\cns\subset \cat{Z}(\cns)$.

Note that $\cns$ is a spherical fusion category, and for any simple $c\in \cns$ we have $[0,2k]c\cong [0,2l]c$ if and only if $k\equiv l \tn{ mod }d$. This means that:
\begin{prop}\label{cnsisvdmod}
$\cns$ is a spherical free $\cat{V}_{d}$-module fusion category.
\end{prop}

\subsubsection{Charged morphisms}

We see from the fusion rules in Equation \eqref{fusionincns} that any object of $\cns$ is a subobject of a tensor product of tensor powers of $[1,1]$ and $[0,2]$. This means that $[1,1]$ generates $\cns$ as a $\cat{V}_d$-module category. Additionally
$$
[1,1]^2=[0,2]\oplus[2,2],
$$
meaning that $[1,1]$ is self-dual up to tensoring by $[0,-2]$, which corresponds to $\underline{\mi 1}\in\cat{V}_d$. To capture this, we use the notion (see eg. \cite{Morrison2019}) of enriching a module category to an enriched category: a category whose hom-sets have been replaced by objects in the monoidal category that acts, see e.g. \cite{Kelly2005,MacLane} for further details. 

We enrich $\cns$ to a $\cat{V}_d$-enriched category $\underleftarrow{\cns}$ by setting, for any $c, c' \in \cns$ and $\underline{k} \in \mathcal{V}_d$:
$$
\underleftarrow{\cns}(c,c')= \bigoplus_{k\in \Z_d} \cns([0,2k] c,c') \underline{k}\in \cat{V}_d.
$$

The inverse to this enrichment is applying the lax monoidal functor $\cat{V}_d(\underline{0},-)$ to the hom-objects. We will use this to pass between $\underleftarrow{\cns}$ and $\cns$ as $\cat{V}_d$-module tensor category freely.

The enrichment satisfies the defining adjunction relation:
$$
\cat{V}_d(\underline{k}, \underleftarrow{\cns}(c,c')) \cong \cns( [0,2k]c,c').
$$

\begin{notation}
We will denote an element of the $\underline{k}$ summand of $\underleftarrow{\cns}(c,c')$ by $f\colon c \rar_k c'$, and the corresponding element of $\cns([0,2k] c,c')$ by $\bar{f}$. 
\end{notation}

Note that $f\colon c \rar_k c'$ and $\bar{f}\colon [0,2k]c\rar c'$ determine each other. We will refer to a morphism $f\colon c \rar_k c'$ as a morphism \textit{of charge $k$}, and to the pairs $f$ and $\bar{f}$ as each other's \textit{mates}. The composition of $f\colon c \rar_k c'$ with $f'\colon c' \rar_{k'} c''$ is $f'\circ f\colon c \rar_{k+k'} c''$, determined by:
$$
\overline{f'\circ f} = \left([0,2(k+k')] c \cong [0,2k'][0,2k]c \xrightarrow{\bar{f}} [0,2k']c' \xrightarrow{\bar{f}'} c''  \right).
$$

The category $\underleftarrow{\cns}$ is monoidal in the sense of \cite[Definition 2.1]{Morrison2019}. On objects the monoidal structure is that of $\cns$. On morphisms, the monoidal product of two morphisms $f_1\colon c_1 \rar_{k_1} c_1'$ and $f_1\colon c_2 \rar_{k_2} c_2'$ is the mate to $f_1\otimes f_2\colon c_1 \otimes c_2 \rar_{k_1+k_2} c_1' \otimes c_2'$,
$$
\xymatrix{
\overline{f_1\otimes f_2} \colon & [0,2(k_1+k_2)]c_1c_2 \cong [0,2k_1][0,2k_2]c_1 c_2 \ar[d]^{\beta_{[0,2k_2],c_1}} & \\  & [0,2k_1]c_1[0,2k_2] c_2 \ar[r]^-{\bar{f}_1 \cdot \bar{f}_2} & c_1' c_2'
}
$$
where $\beta$ is the braiding in $\cns$.

\subsubsection{Charged duality}

In this language, $[1,1]$ is self-dual of charge $-1$. That is, we have an evaluation $\ev\colon [1,1]^2 \rar_{-1} [0,0]$ exposing $[1,1]$ as its own left dual. It is defined by
$$
\overline{ev}\colon [0,-2][1,1][1,1] \cong [1,-1][1,1] \xrightarrow{e} [0,0],
$$
where $e$ is an evaluation map witnessing $[1,-1]$ as the left dual of $[1,1]$. The co-evaluation has charge $1$ and is the morphism $\coev\colon [0,0] \rar_{1} [1,1]^2$ defined through
$$
\overline{coev}\colon [0,2][0,0] \xrightarrow{\id_{[0,2]} \cdot ce} [0,2][1,1][1,-1] \xrightarrow{\beta_{[0,2],[1,1]}} [1,1][0,2][1,-1] \cong [1,1]^2,
$$
where $ce$ is the co-evaluation associated to $e$, and $\beta_{[0,2],[1,1]}$ is the braiding. The isomorphism $[0,2][1,-1] \cong [1,1]$ is chosen to be the dual element under the composition pairing to the isomorphism $ [1,1] \cong [0,2][1,-1]$ used in the definition of the evaluation. A quick computation shows that these morphisms satisfy the usual snake relations. For the right duality, we have $\widetilde{ev}\colon [1,1]^2 \rar_{-1} [0,0]$ defined by
$$
\overline{\widetilde{\ev}}\colon [0,-2][1,1][1,1] \xrightarrow{\beta_{[0,-2],[1,1]}.\id_{[1,1]}} [1,1][0,-2][1,1] \cong [1,1][1,-1] \xrightarrow{p(e)} [0,0],
$$
where $p(e)$ denotes $e$ composed with the pivotal structure. The corresponding co-evaluation is defined by
$$
\overline{\widetilde{\coev}}\colon [0,2][0,0] \xrightarrow{\id_{[0,2]} \cdot p(ce)} [0,2][1,-1][1,1] \cong [1,1]^2.
$$
One computes that $\ev \circ \widetilde{\coev}$ and $\widetilde{\ev}\circ \coev$ agree and are equal to
$$
\qdim([1,1])= 2 \cos(\frac{\pi}{d}),
$$
showing that this number also deserves to be called the quantum dimension of $[1,1]$ in this context.

\subsubsection{Generating $\cns$}

We can now prove:

\begin{thm}\label{cnsgenerated}
The category $\cns$ is equivalent as $\cat{V}_d$-module tensor category to the pivotal free $\cat{V}_d$-module fusion category generated by a single object that is self-dual of charge $-1$ and quantum dimension $2 \cos(\frac{\pi}{d})$.
\end{thm}

\begin{proof}
We need to show that the data of a self-dual object of charge $-1$ and quantum dimension $\kappa= 2 \cos (\frac{\pi}{2})$ is sufficient to recover $\cns$ as a $\cat{V}_d$-module tensor category. To do this, we will generate a spherical free $\cat{V}_d$-module tensor category $\cat{K}$ from this data, and provide an equivalence between the maximal fusion quotient of its idempotent completion and $\cns$.

We will closely follow the classical treatment of Temperley-Lieb categories as outlined in Section \ref{temperleysection}. The difference is that we will generate a $\cat{V}_d$-enriched and tensored category to accommodate the charge. The objects of $\cat{K}$ are pairs $(k,n)\in \Z_d \times \mathbb{N}$ (which will eventually be mapped to $[0,2k]\otimes [1,1]^{\otimes n}$). On objects, the tensor product on $\cat{K}$ is just: 
$$
(k,n) \otimes (l,m)= (k+l,n+m).
$$
The morphisms in this category are generated from the duality data $\ev$ and $\coev$ and by enforcing that the resulting category is a free $\cat{V}_d$-module category. That is, morphisms from $(k,n)$ to $(l,m)$ are a tensor product of $\underline{l-k}$ with the span of Kauffman diagrams $K_{m,n}$ between $m$ and $n$ up to isotopy:
$$
\cat{K}((k,n),(l,m))= \underline{l-k}\otimes\C[K_{m,n}]
$$
Here, we make $\C[K_{m,n}]$ into a $\Z_d$-graded vector space by setting the degree of a diagram to be the number of co-evaluations minus the number of evaluations in the diagram. This encodes the charge of the duality data. Observe that $\C[K_{m,n}]$ is trivial when $m-n$ is odd, and has degree $\frac{m-n}{2}$ otherwise. 

Composition of morphisms from $(k,n)$ to $(l,m)$ with morphisms from $(l,m)$ to $(j,r)$ is induced from the map $K_{m,n}\times K_{r,m}\rar K_{r,n}$ that glues diagrams, together with the obvious isomorphism $\underline{l-k}\otimes \underline{r-l}\cong \underline{r-k}$. The tensor product of morphisms is induced from juxtaposition of diagrams.

We now want to compute the idempotent completion of $\cat{K}$. To do this, we need to examine the endomorphism algebras in $\cat{K}$. The endomorpisms of $(k,n)$ are of charge 0, meaning that finding all the idempotents in $\cat{K}$ reduces to the classical (uncharged) case. The hom-object $\cat{K}((k,n),(k,n))$ is for each $k$ the classical Temperley-Lieb algebra $TL_n(\kappa)$. Finding the idempotents then proceeds like in Section \ref{temperleysection}. 

This gives for each $k$ a sequence $\{p_{k,n}\in \End((k,n))\}$, giving rise to objects $\langle\!\langle k,n \rangle\!\rangle\subset (k,n)$ with $0\leq n \leq d-2$ in the idempotent completion of $\cat{K}$. 

Observe that for any $k$ and $l$ we have $(k,n)\cong (l,n)$ along an isomorphism of charge $l-k$. These isomorphisms descend to isomorphisms
\begin{equation}\label{chargeisos}
\langle\!\langle k,n \rangle\!\rangle \xrightarrow{\cong}_{l-k} \langle\!\langle l, n \rangle\!\rangle, 
\end{equation}
giving for each $n$ the subcategory spanned by the $\langle\!\langle k,n \rangle\!\rangle$ the structure of a free rank one module category over $\cat{V}_d$.

In the classical charge 0 case, one uses the Kauffman diagrams to find the complementary summands to $\langle\!\langle k,n \rangle\!\rangle$ in $ (k,n)$. In our charge $-1$ case, this analysis stays almost the same\footnote{This should not come as a surprise. Idempotent completion only sees the underlying linear category, and does not heed the charges.}, the only difference being that the inclusions can now carry charge. Taking this into account and using Equation \eqref{chargeisos} to obtain a morphism of charge zero, we find
$$
\langle\!\langle k,n \rangle\!\rangle \otimes \langle\!\langle l,1 \rangle\!\rangle \cong \langle\!\langle k+l-1,n-1 \rangle\!\rangle \oplus \langle\!\langle k+l, n+1 \rangle\!\rangle,
$$
where the projection onto the first summand is constructed from $p_{n-1}$ and $\ev$ together with a charge $1$ isomorphism.

Just like in the classical case \cite[Proposition 2.1]{Goodman2003} the final well-defined projectors $p_{k,d-1}$ generate the tensor ideal of negligible morphisms. This means we just need to provide a self-dual object in $\cns$ of charge $-1$, quantum dimension $\kappa$ and vanishing Jones-Wenzl projector. By the discussion preceding the Theorem, $[1,1]\in \cns$ is such an object.
\end{proof}

\begin{rmk}\label{drcrcompareremark1}
We worked with the object $[1,1]$ for concreteness, but any $[1,m]$ with $m$ odd generates $\cns$ as a $\cat{V}_d$-module category, and is self dual of charge $-m$ with quantum dimension $2\cos(\frac{\pi}{d})$. This gives a presentation of $\cns$ for each odd number in between $0$ and $2d$. If $d$ is itself odd, this means that the object $[1,d]$ is self-dual of charge $d\equiv 0 \tn{ mod }d$ and generates $\cns$. Running through the proof of Theorem \ref{cnsgenerated} for this object, we recover the observation from \cite[Proposition 2.1]{Davydov2018} that $\cns$ is, for odd $d$, the Deligne product of the fusion quotient of a Temperley-Lieb category with $\cat{V}_d$.
\end{rmk}

\section{The Landau-Ginzburg side}\label{Section3}
We now turn our attention to the Landau-Ginzburg side of the correspondence. This side is described using matrix factorisations, and we briefly recall some details about these. After that, we define $\pd$, and prove that it can be given the structure of a $\vd$-module tensor category.

\subsection{Matrix factorisations}
Here we introduce some basics of matrix factorisations. Denote by $R$ a commutative $\mathbb{C}$-algebra, and let us pick $W \in R$.

\begin{df}
A \textit{matrix factorisation} of $W$ consists of a pair $\left( M,d^M \right)$, where:
\begin{itemize}
    \item[-] $M$ is a finite-rank, free, $\mathbb{Z}_2$-graded $R$-module $M=M_0 \oplus M_1$; and
    \item[-] $d^M \colon M \to M$ is a degree 1 $R$-linear endomorphism $d^M=\left( \begin{matrix} 0 & d_1^M \\ d_0^M & 0 \end{matrix} \right)$ called the \textit{twisted differential}, satisfying $d^M \circ d^M=W.\mathrm{id}_M$.
\end{itemize}
\end{df}

If clear from the context, we will simply use $M$ to denote a matrix factorisation $\left( M,d^M \right)$. If $(R',W')$ is another ring with chosen element, a matrix factorisation for $R\otimes_{\mathbb{C}} R'$ and potential $W\otimes 1-1\otimes W'$ is called a \emph{matrix bifactorisation} following \cite{Carqueville2010}. We will be interested in the case where $R=R'$ and $W=W'$: bifactorisations of $(R,W)$. The underlying module $M$ for such a bifactorisation of $W$ is an $R-R$-bimodule, and $d^M \circ d^M=W.\mathrm{id}_M-\mathrm{id}_M.W$. Note here that a $R$-bimodule is free if it is free as an $R \otimes_{\mathbb{C}} R$-module.

\begin{df}
Let $\left( M,d^M \right)$, $\left( N,d^N \right)$ be two matrix bifactorisations of $W$. We define a \textit{morphism of matrix factorisations} to be an $R$-bilinear morphism $f \colon M \to N$. Note that these morphisms inherit the $\mathbb{Z}_2$-grading of the base module, so we can write them in components as $f=\left( f_0,f_1 \right)$.
\end{df}

The category $\mathrm{MF}_{\mathrm{bi}} \left( W \right)$ with matrix bifactorisations as objects and morphisms as defined above has the structure of a differential $\mathbb{Z}_2$-graded category \cite{Weibel1995}, with differential in the morphism space defined as follows: let $f \in \mathrm{Hom}_{\mathrm{MF_{\mathrm{bi}}(W)}} \left( M,N \right)$, then
$$\delta \left( f \right):=d^N \circ f -\left( -1 \right)^{|f|} f \circ d^M.$$

We will work in the associated homotopy category:
\begin{df}\label{hmfdef}
The category $\mathrm{HMF}_{\mathrm{bi}} \left( W \right)$ is the category with the same objects as $\mathrm{MF}_{\mathrm{bi}} \left( W \right)$, and as hom-objects the zero degree homology of the hom-objects of $\mathrm{MF}_{\mathrm{bi}} \left( W \right)$. 
\end{df}

Assuming that $R=\mathbb{C} \left[ x_1,\ldots,x_n \right]$ and $W$ is a polynomial satisfying that 
$$
\mathrm{dim}_\mathbb{C} \left( \frac{R}{\langle \partial_{x_1} W,\ldots,\partial_{x_n} W \rangle} \right)
$$
is finite (i.e. a \textit{potential}) the category $\mathrm{HMF}_{\mathrm{bi} \left( W \right)}$ is tensor \cite{Carqueville2016c,Carqueville2010}:

\begin{df}
Given $\left( M,d^M \right)$, $\left( N,d^N \right) \in  \mathrm{HMF}_{\mathrm{bi}} \left( W \right)$, their \emph{tensor product} $M\otimes N$ is defined as follows: the base module is $M \otimes_R N$ and the twisted differential is $d^{M \otimes N}=d^M \otimes \mathrm{id}_N+\mathrm{id}_M \otimes d^N$. 
\end{df}
This tensor product defines a monoidal structure with trivial associators, we will discuss what is needed about the unitors below. In terms of the $\mathbb{Z}_2$-grading, the underlying module is 
$$M \otimes_R N=\left( M_0 \otimes_R N_0 \oplus M_1 \otimes_R N_1 \right) \oplus \left( M_1 \otimes_R N_0 \oplus M_0 \otimes_R N_1\right).
$$ 
To see that the differential on the tensor product squares to the desired potential, one needs to use that the composition of tensor products of graded morphisms follows the Koszul sign rule. 

\begin{notation}\label{gradingcomponentsnotation}
In this paper, we will often need to consider a morphism $f$ in $\mathrm{HMF}_{\mathrm{bi} \left( W \right)}$ between a matrix factorisation $M'$ and a tensor product $M\otimes N$ of two matrix factorisations. Both when $M\otimes N$ is the source and when it is the target such a morphism has four components: the maps $M_{i}\otimes N_j\leftrightarrow M'_{i+j}$ with $i,j\in \Z_2$. In the rest of this paper, we will denote these components by $ f_{ij}$.
\end{notation}

For further details on the (higher) categorical structure of $\mathrm{HMF}_\mathrm{bi} \left( W \right)$, we refer to \cite{Carqueville2016c,Carqueville2010,Carqueville2012}.

\subsection{The category $\mathbf{\mathcal{P}_d}$}
Let us fix from now on $R=\mathbb{C} \left[ x \right]$ and $W=x^d$ for $d \in \mathbb{Z}_{\geq 2}$. We are interested in a particular subcategory $\mathbf{\mathcal{P}_d}$ of $\mathrm{HMF}_{\mathrm{bi}} \left( W \right)$ consisting of so-called permutation type matrix factorisations. After explaining what these are we examine the structure of this category as fusion category and prove that it admits a $\vd$-module tensor category structure.

\subsubsection{Permutation type matrix factorisations}
Consider the matrix bifactorisation in $\mathrm{HMF}_{\mathrm{bi}} \left( W \right)$ with base module $P=P_0 \oplus P_1=\mathbb{C} \left[ x,y \right]^{\oplus 2}$ and twisted differential given by the matrix 
$$
d_S=\left( \begin{matrix} 0 & \prod\limits_{j \in S} \left( x-\eta^j y \right) \\ \frac{x^d-y^d}{\prod\limits_{j \in S} \left( x-\eta^j y \right)} & 0 \end{matrix} \right),
$$ 
where $S \subset \lbrace 0,\ldots,d-1 \rbrace$. Here, $\eta$ can be chosen to be any primitive $d$-th root of unity, for us this is $\eta=e^{\frac{2\pi i}{d}}$, cf. Section \ref{notationsection}. We will denote such a matrix factorisation as $P_S$ and we will call it a \textit{permutation type matrix factorisation}. Note that the tensor unit of $\mathrm{HMF}_{\mathrm{bi}} \left( W \right)$ is the permutation type matrix factorisation $I=P_{\lbrace 0 \rbrace}$. 

\begin{notation}
We set $S=\lbrace m,\ldots,m+l \rbrace$ and denote the associated permutation type matrix factorisations as $P_{m;l}\equiv P_S$. We will use $P_S$ and $P_{m;l}$, and $S$ and $m;l$ interchangeably.
\end{notation}

In what follows, we will only be interested in the full subcategory $\mathcal{P}_d$ of $\mathrm{HMF}_{\mathrm{bi}} \left( W \right)$ spanned by the simples $P_{m;l}$. This is a semi-simple category \cite{Davydov2018}. Many of the results listed here are specialised to $\mathcal{P}_{d}$ from results for $\mathrm{HMF}_{\mathrm{bi}} \left( W \right)$ for any $W$. For the general theory see \cite{Carqueville2010}.

\subsubsection{The monoidal structure of $\pd$}
The subcategory $\mathcal{P}_{d}$ is closed under tensor product (see \cite{Brunner2007,Davydov2018}). The fusion rules are 
\begin{equation}\label{fusionrulesPd}
P_{m;l} \otimes P_{m';l'} \cong \bigoplus\limits_{\nu=\vert l-l' \vert \mathrm{ step } 2}^{\mathrm{min\left( l+l',2d-4-l-l' \right)}} P_{m+m'+\frac{1}{2} \left( l+l'-\nu \right);\nu}    
\end{equation}
for $l,l' \in \lbrace 0,\ldots,d-2 \rbrace$, and $m,m' \in \mathbb{Z}_{2d}$. The objects $P_{m;0}$ span a copy of $\Vect[\mathbb{Z}_d]$ with trivial associator, this is the tensor category underlying $\cat{V}_d$ (see Lemma \ref{jscomputation}). This gives an embedding of tensor categories
\begin{equation}\label{vdintopd}
    \vd \hookrightarrow \pd,
\end{equation}
sending the one dimensional degree-$m$ vector space $\C_m\in \vd $ to $\underline{m}:=P_{m;0}$. We will eventually lift this embedding to a central functor $\vd\rar \cat{Z}(\pd)$, this is Theorem \ref{pdisvdmod}.

The left and right unitors $\lambda_S:=\lambda_{P_S} \colon P_{0;0} \otimes P_S \to P_S$ and $\rho_S:=\rho_{P_S} \colon P_S \otimes P_{0;0} \to P_S$ are given by (recall Notation \ref{gradingcomponentsnotation})
\begin{equation}
   (\lambda_{S})_{ij} := \delta_{i,0}L_{(P_S)_i} \mbox{ and }
    (\rho_{S})_{ij} := \delta_{j,0} R_{(P_S)_j},
    \label{unitors}
\end{equation}
Here the maps $L_P$ and $R_P$ are defined for a $\mathbb{C}[x]$-bimodule $P$ as
\begin{align*}
L_P \colon \mathbb{C} \left[x,y \right] \otimes P & \to P \\ f \left( x,y \right) \otimes p & \mapsto f \left( x,x \right).p \\ R_P \colon P \otimes  \mathbb{C} \left[ x,y \right] &\to P \\ p \otimes f \left( x,y \right) & \mapsto p. f \left( x,x \right).
\end{align*}

\subsubsection{Duals in $\pd$}

Any finite-rank matrix factorisation $M$ in $\mathrm{HMF}_{\mathrm{bi}} \left( W \right)$ has a left dual $M^\dagger$ (see \cite{Carqueville2016c,Carqueville2012}) with the same underlying module $M$ and twisted differential $d^{M^\dagger}=-d^M$. For a permutation type matrix factorisation $P_S\in \mathcal{P}_d$, there is an isomorphism
\begin{equation*}
    \left( P_S \right)^\dagger \simeq P_{-S} 
\end{equation*}
given by multiplication by $(-1)^{|S|+1}\prod\limits_{j\in S}\eta^{j}$ on the odd summand and the identity on the even summand \cite[Equation 3.4]{Davydov2018}. 
Combining this isomorphism with the evaluation and coevaluation maps from \cite[Page 12]{Davydov2018} gives maps 
\begin{align*}
    \mathrm{coev}_S:=\mathrm{coev}_{P_{S}}  \colon& I \to  P_{S}  \otimes P_{-S} \\
    \mathrm{ev}_{S}:=\mathrm{ev}_{P_S}  \colon&  P_{-S} \otimes P_{S}  \to I 
\end{align*}
exhibiting $P_{-S}$ as the left dual $P_S^*$ of $P_S$. This defines a functor $(-)^*\colon \pd^\tn{op}\rar \pd$. The explicit form of the coevaluation map is in components (see Notation \ref{gradingcomponentsnotation}):
\begin{equation}
       (\mathrm{coev}_{S})_{ij}=\begin{cases}(-1)^{|S|+1}\prod\limits_{k\in S}\eta^{k}& \mbox{ for }i=0, j=1 \\
                                            1 & \mbox{ for }i=1,j=0\\
                                        \left((-1)^{|S|+1}\prod\limits_{k\in S}\eta^{k}\right)^{j}\frac{d_{i+1}^{S}(x,y)-d_{i+1}^{S}(z,y)}{x-z}& \mbox{ for } i+j=0
    \end{cases} 
     \label{coevforpermtype}
\end{equation}
The evaluation map has components that are maps $\C[x,y,z]\rar \C[x,y]$ given by 
\begin{equation}\label{evforpermtype}
(\ev_{S})_{ij}(f(x,y,z))= \begin{cases}- \mathcal{G}_S(f(x,y,z)& \mbox{ for }i=j=0\\
                                    0 & \mbox{ for }i=j=1\\
                                    (-1)^{|S|+1}\prod\limits_{k\in S}\eta^{-k}\frac{\mathcal{G}_S(f(x,y,z)}{x-z}& \mbox{ for }i=1,j=0\\
                                    -f(x,0,z) &\mbox{ for }i=0,j=1 \end{cases}.    
\end{equation}
The map $\mathcal{G}_S$ is defined using a contour integral over $y$ and is $\C[x,z]$-linear, see \cite{Davydov2018}. For our purposes it suffices to know that 
\begin{equation}\label{gmapvalues}
\mathcal{G}_{m;0}(1)=-1, \quad \mathcal{G}_{m;1}(1)=0 \tn{ and }\quad \mathcal{G}_{m;1}(y)=1,    
\end{equation}
as one can see from a brief computation. For further details we refer to \cite{Carqueville2012,Davydov2018}. Note that $\pd$ is a semi-simple rigid tensor category with simple unit: a fusion category \cite{Etingof2002}.

\subsubsection{A spherical structure for $\mathcal{P}_d$}
Choosing the left dual $P_S^*$ to a simple $P_S\in \mathcal{P}_d$ to be $P_{-S}$ induces a pivotal structure by identifying $P_S^{**}=(P_{-S})^*=P_{S}$ along the identity map. With this pivotal structure we can compute the quantum dimensions of simple objects $P_S$ by computing the composite 
$$
\qdim_{L}(P_S)=ev_S\circ\coev_{-S}.
$$ 
We only need to evaluate either the even or the odd part of this map. Furthermore, note (Equation \eqref{evforpermtype}) that $(\ev_{S})_{11}=0$. So it suffices to compute (using that $\mathcal{G}_{S}$ is $\C[x,z]$-linear and Equation \eqref{coevforpermtype})
$$
(\ev_{S})_{00}\circ(\coev_{-S})_{00}=-\mathcal{G}_{S}(\frac{d^{\mi S}_1(x,y)-d^{\mi S}_1(z,y)}{x-z}).
$$
We will only need
\begin{align*}
    \qdim_L(P_{m;0})&=-\mathcal{G}_{m;0}(\frac{x-z}{x-z})=1\\
    \qdim_L(P_{m;1})&=-\mathcal{G}_{m;1}(\frac{x^2+(\eta^{-m}+\eta^{-m-1})(z-x)y-z^2}{x-z})=\eta^{-m}+\eta^{-m-1},
\end{align*}
where we used Equation \eqref{gmapvalues}. As $\qdim_R(P_{m;l})$ is the complex conjugate of $\qdim_L(P_{m;l})$, this pivotal structure is not spherical. This is undesirable: $\cns$ is spherical. Note however that $\pd$ is a pseudo-unitary fusion category (its fusion matrices are those of $\repsu$ extended by zeros), and therefore admits a unique spherical structure \cite[Proposition 9.5.1]{Etingof2015}. By \cite[Proposition 2.4]{Muger2003} this spherical structure is such that
$$
\qdim(P_{m;l})^2=\qdim_L(P_{m;l})\qdim_R(P_{m;l}),
$$
where $\qdim$ denotes the dimension in the spherical structure. In particular we have that
\begin{equation}
    \label{quantumdimensionofP01}
    \qdim(P_{m;1})=2\cos (\frac{\pi}{d})
\end{equation}
in this spherical structure. Henceforth we will use this spherical structure.

\subsection{Fusion isomorphisms for $\mathcal{P}_d$}
In order to prove that $\pd$ admits the structure of a $\vd$-module tensor category, we will need concrete choices for the isomorphisms in Equation \eqref{fusionrulesPd}. We further need to know how these choices interact under composition. 

There is a $\Z_d$-action on $\pd$ that helps choosing the fusion isomorphisms. We will refer to this action as twisting, and describe it in the next section.

\subsubsection{Twisting and untwisting matrix factorisations}
The algebra automorphism $\sigma \colon \mathbb{C}\left[ x \right] \to \mathbb{C}\left[ x \right]$ defined by $x \mapsto \eta x$ leaves the potential $W=x^d$ invariant. This induces a group isomorphism $\mathbb{Z}_d \cong \mathrm{Aut} \left( \mathbb{C} \left[ x \right], x^d \right) $, given by $k \mapsto \sigma^k$. These automorphisms then act on the bifactorisations of $x^d$, we write $_a \left( P_S \right)_b$ to denote a permutation type matrix factorisation whose underlying $\mathbb{C} \left[ x \right]$-bimodule is equal to its original one as a $\mathbb{Z}_2$-graded $\mathbb{C}$-graded vector space, but has its left and right actions twisted by $\sigma$ in the following way. For $p \in \mathbb{C} \left[ x \right]$, and $m \in \mathbb{C} \left[ x,y \right]$ we send
\begin{align*}
\left( p,m \right) & \mapsto \sigma^{-a} \left( p \right).m , \\
\left( m,p \right) & \mapsto m.\sigma^{b} \left( p \right) .
\end{align*}
where . denotes the left or right action on the original bimodule. 

It will often be convenient to keep track of which variable we are acting on by setting
\begin{align*}
\sigma_{x_i}^a\colon &\C[x_1,\dots, x_n] \rar \C[x_1,\dots ,x_n]\\
&f(x_1, \dots, x_n) \mapsto f(x_1, \dots , \eta^a x_i, \dots, x_n).    
\end{align*}

Twisting $P_S$ can be untwisted to give $P_{S'}$ with $S'$ a $\Z_d$-translate of $S$:

\begin{df}
$s_{a,b} \colon P_{S-a-b} \to _a \left( P_S \right)_b$ is the isomorphism of matrix factorisations given in components as:
$$s_{a,b}:=\left( \sigma^{-a} \otimes \sigma^{b}, \eta^{- \vert S \vert a} \sigma^{-a} \otimes \sigma^{b}\right).$$
\end{df}

We will need a couple of facts about this twisting and untwisting.
\begin{facts}\label{mainfacts}
\begin{enumerate}
    \item (See Lemma 3.5, \cite{Davydov2018}) Note that $\forall a,b \in \mathbb{Z}_d$, $\,_{a}\left( - \right)_b$ defines an autoequivalence of $\mathcal{P}_d$. If $b=-a$, this auto-equivalence is tensor. When either $a$ or $b$ is zero, then we will abbreviate the notation by simply adding a subscript at the non-zero side, e.g. $\,_{a}\left( \mathbbm{1} \right)_0=\,_{a}\mathbbm{1}$. Note that we can relate any $\underline{a}$($=P_{a;0}$) and $\,_{\mi a}\mathbbm{1}$ via $s$: we have $s_{\mi a,0} \left( \underline{a} \right)=\,_{\mi a}\mathbbm{1}$.
    \item\label{fact2} A direct consequence of this which we will use later is that $\underline{(a+b)}$ is isomorphic to $\,_{\mi a}\mathbb{I}_{\mi b}$ along $s_{\mi a, \mi b}= \,_{\mi a}(s_{0,\mi b})\circ s_{\mi a, 0}$, and similarly $\,_{\mi a \mi b}\left( P_{S} \right)\cong \,_{\mi a}\left(P_{S+b}\right)$ along $\,_{\mi a}(s_{\mi b;0})^{-1}= s_{\mi a,0}\circ s^{-1}_{\mi a \mi b,0}$. 
    \item Note that for any two underlying modules of two permutation type matrix factorisations respectively $M$ and $N$, twisting the intermediate variable can be transferred through the tensor product:
$$
M_a\otimes N = M\otimes_a N.
$$
    \item The maps $s$ respect the addition in $\Z_d$ in the sense that
\begin{equation}\label{sisadditive}
    s_{a+b,-a-b}=s_{b,-b}s_{a,-a}=s_{a,-a}s_{b,-b}.
\end{equation}
\end{enumerate}
\end{facts}

In what follows we will need the following relationship between the unitors and $s$:

\begin{lem}\label{sandunitorinteraction}
For $a\in \Z_d$ and any $P_{S} \in  \mathcal{P}_d$, we have:
\begin{align*}
\left(\,_{a}\mathbb{I} P_{S} \xrightarrow{\,_a( s_{\mi a,a} ) \cdot \id_{P_{S}}} \mathbb{I}_a P_{S} \xrightarrow{\lambda_{\,_a(P_{S})} } \,_a P_{S} \right) &= \left(\,_a\mathbb{I}  P_{S} \xrightarrow{\,_a(\lambda_{P_{S}})}\,_aP_{S}\right)\\
\left( P_{S}\mathbb{I}_a \xrightarrow{ \id_{P_{S}}\cdot ( s_{\mi a,a} )_a }  (P_{S})_a\mathbb{I} \xrightarrow{\rho_{(P_{S})_a} }  (P_{S})_a \right) &= \left( P_{S}\mathbb{I}_a \xrightarrow{(\rho_{P_{S}})_a}(P_{S})_a\right).
\end{align*}
\end{lem}
\begin{proof}
Note that for the degree $i$ summand of $P_{S}$ the left hand side acts on $f(x,y)\otimes p \in \C[x,y] \otimes_{\C[y]} P_i$ as
$$
f(x,y)\otimes p \xmapsto{\sigma_x^a \cdot \sigma_y^a \cdot \id_{p}} f(\eta^a x, \eta^a y)\otimes p \xmapsto{L_{P_i}} f(x,x).p, 
$$
where $.$ denotes the untwisted action of $\C[x]$ on the $\C[x]$-bimodule $P_i$, and we have used that the twisted action $._a$ of $\C[x]$ on $\,_aP_i$ is given by $x._a(-)= \eta^{\mi a}x.(-)$. This agrees with the definition of the unitor from Equation \eqref{unitors}. The proof of the second equality is analogous. 
\end{proof}

\subsubsection{Explicit fusion isomorphisms}\label{fusionsection}
The fusion rules from Equation \eqref{fusionrulesPd} are established in \cite[Appendix B]{Davydov2018} by finding explicit morphisms exhibiting the direct sum decompostion of two simples into summands. This gives a basis for the hom-spaces  $\mathcal{P}_d(P_{S}P_{S'},P_{S''})$ and the goal of this section is to establish how this interacts with the associators. The associators in $\mathcal{P}_d$ are trivial, meaning that they induce the identity map between  $\mathcal{P}_d((P_{S}P_{S'})P_{S''},P_{S'''})$ and  $\mathcal{P}_d(P_{S}(P_{S'}P_{S''}),P_{S'''})$. Our aim is to examine what this identity map looks like in the bases induced by composition of the bases from \cite[Appendix B]{Davydov2018}. This sets us up to give a coherent set of isomorphisms $\underline{a}P_{m;l}\underline{-a}\cong P_{m;l}$ in Section \ref{pdasmodcatsection}.

Establishing the results we need requires examining the underlying modules and the maps between them in detail, leading to a couple of lengthy arguments. Fortunately, the results can be summarised in simple string diagrams, and we work with those in the rest of the paper. 

Following \cite[Appendix B]{Davydov2018} we take
$$
\psi^L_{a;S}\colon \underline{a}P_S \xrightarrow{ s_{0;\mi a} \cdot \mathrm{id}_{P_S}} \mathbb{I}_{\mi a}P_S\xrightarrow{\id_{\mathbb{I}}\cdot(s_{\mi a,0})^{-1}} \mathbb{I}P_{S+a} \xrightarrow{\lambda_{P_{S+a}}} P_{S+a}
$$
(also denoted as $\psi^L_{a,(m;l)}$ if we choose to make the permutation set explicit) as the basis for $\mathcal{P}_d (\underline{a}P_{S},P_{S+a})$. By naturality of the unitors, we can swap the last two morphisms in this composite. In string diagrams, we will depict this morphism as:
$$
\psi^L_{a;S}=
\hbox{
\begin{tikzpicture}[baseline=(current  bounding  box.center)]
\node (m) at (-0.75,0){$\underline{a}$};
\node (p) at (0,0){$P_{S}$};
\coordinate (v) at (0,0.7);
\node (pm) at (0,1.3){$P_{S+a}$};

\draw[thick,blue] (m) to [out=north,in=west] (v);
\draw[thick] (p) to (v) to (pm);
\draw[fill,blue] (v) circle [radius=0.06];
\end{tikzpicture}
}.
$$

To aid readability, we will use blue strands to indicate objects in $\cat{V}_d\subset \mathcal{P}_d$. This choice of basis behaves well under composition:

\begin{lem}\label{leftpsiasso}
Let $a,b\in \Z_d$ and $P_{m;l}\in \mathcal{P}_d$. Then
$$
\psi^L_{a;S+b}\circ(\id_{\underline{a}} \otimes \psi^L_{b;S}) = \psi^L_{a+b;S}\circ(\psi^L_{a;(b;0)}\otimes\id_{P_{S}}),
$$
or in string diagrams:
$$
\hbox{
\begin{tikzpicture}[baseline=(current  bounding  box.center)]

\node (a) at (-0.8,0){$\underline{a}$};
\node (b) at (-0.5,0){$\underline{b}$};
\node (p) at (0,0){$P_{S}$};
\coordinate (v) at (0,0.7);
\coordinate (v1) at (0,1);
\node (pm) at (0,1.6){$P_{S+a+b}$};

\draw[thick,blue] (a) to [out=north,in=west] (v1);
\draw[thick,blue] (b) to [out=north,in=west] (v);
\draw[thick] (p) to (v) to (v1) to (pm);
\draw[fill,blue] (v) circle [radius=0.06];
\draw[fill,blue] (v1) circle [radius=0.06];
\end{tikzpicture}
}
=
\hbox{
\begin{tikzpicture}[baseline=(current  bounding  box.center)]

\node (a) at (-0.8,0){$\underline{a}$};
\node (b) at (-0.5,0){$\underline{b}$};
\node (p) at (0,0){$P_{S}$};
\coordinate (v) at (-0.5,0.6);
\coordinate (v1) at (0,1);
\node (pm) at (0,1.6){$P_{S+a+b}$};

\draw[thick,blue] (a) to [out=north,in=west] (v);
\draw[thick,blue] (b) to [out=north,in=south] (v) to [out=north, in=west] (v1);
\draw[thick] (p) to (v1) to (pm);
\draw[fill,blue] (v) circle [radius=0.06];
\draw[fill,blue] (v1) circle [radius=0.06];
\end{tikzpicture}
}.
$$

\end{lem}
\begin{proof}
Using the interchange law for the tensor product and naturality of the unitors, we can re-arrange both sides so that they start with 
$$
s_{0,\mi a}\otimes \id_{\underline{b} P_{S}} \colon \underline{a}\underline{b}P_{S} \rar \mathbb{I}_{\mi a} \underline{b}P_{S}
$$
and end with $\lambda_{P_{S+a+b}}$. As both these maps are isomorphisms, we can cancel them, leaving us with

$$
\mathbb{I}_{\mi a}\underline{b} P_{S} \xrightarrow{\id_{\mathbb{I}_{\mi a}} \cdot s_{0,\mi b} \cdot \id_{P_S}} \mathbb{I}_{\mi a}\mathbb{I}_{\mi b}P_{S} \xrightarrow{\id_{\mathbb{I}_{\mi a}} \cdot s^{-1}_{\mi b, 0}}  \mathbb{I}_{\mi a}\mathbb{I}P_{S+b} \xrightarrow{\id_{\mathbb{I}_{\mi a}} \cdot \lambda_{P_{S+b}}}\,_{\mi a}\mathbb{I}P_{S+b}\xrightarrow{s^{-1}_{\mi a,0}} \mathbb{I}P_{S+a+b} 
$$

for the left hand side, and

$$
\mathbb{I}_{\mi a}\underline{b} P_{S}\xrightarrow{s^{-1}_{\mi a,0} \cdot \id_{P_S}} \mathbb{I}\underline{(a+b)}P_{S}\xrightarrow{\id_{\mathbb{I}} \cdot s_{0,\mi a\mi b} \cdot \id_{P_S}} \mathbb{I} \mathbb{I}_{\mi a \mi b} P_{S}\xrightarrow{\lambda_{\mathbb{I}_{\mi a \mi b} P_{S}}}\mathbb{I}_{\mi a \mi b}P_{S} \xrightarrow{s^{-1}_{\mi a \mi b,0}}\mathbb{I}P_{S+a+b} 
$$

for the right hand side. Both sides are an image under the invertible functor $\mathbb{I}\otimes -$, so we can cancel the leftmost $\mathbb{I}$. Fact \ref{mainfacts}.\ref{fact2} allows us to form commutative squares on the outermost two terms in both composites. This leaves to prove that
\begin{center}
    \begin{tikzcd}[column sep=large]
        \,_{\mi a}\mathbb{I}_{\mi b}P_{S} \arrow[r, "\,_{\mi a}\lambda"]\arrow[d,"s^{-1}_{\mi a, \mi b} \cdot \id_{P_S}"]     & \,_{\mi a\mi b} \left( P_{S} \right) \arrow[r,"\,_a(s^{-1}_{\mi b, 0})"]\arrow[d,"\lambda_{\,_{\mi a\mi b} \left( P_{S} \right) }^{-1}"] & \,_a\left(P_{S+b} \right) \arrow[d,"\,_{\mi a}(s_{\mi b;0})"] \\
        \underline{(a+b)}P_{S}  \arrow[r,"s_{0,\mi a \mi b}"'] & \mathbb{I}_{\mi a \mi b}P_{S} \arrow[r, "\lambda_{\,_{\mi a\mi b} \left( P_{S} \right) }"'] & \,_{\mi a \mi b } \left(P_{S}\right)
    \end{tikzcd}
\end{center}
commutes, where we used naturality of $\lambda$ to swap the middle maps in the top row, and inserted the inverse unitor. The rightmost square in the diagram clearly commutes, so we just need to show the leftmost square commutes. Using that $s_{0;\mi a \mi b}s_{\mi a, \mi b}^{-1}= \,_{\mi a}(s_{a,\mi a})$, we see this comes down to showing $\lambda\circ \,_{\mi a}(s_{a,\mi a})= \,_{\mi a}\lambda $. This is Lemma \ref{sandunitorinteraction}, so we are done.
\end{proof}

This implies that the associators are given by the identity matrix in this basis. In the definition of $\psi^L_{a;S}$ we used the left unitor. There is a mirrored version
$$
\psi^R_{S;a} \colon P_{S} \underline{a} \xrightarrow{\mathrm{id}_{P_S} \cdot s_{-a,0}} \left(P_{S}\right)_{\mi a} \mathbb{I} \xrightarrow{s^{-1}_{0,\mi a}\cdot \mathrm{id}_{\mathbb{I}}}  \left(P_{S+a}\right) \mathbb{I} \xrightarrow{\rho_{\left(P_{S+a}\right)}} P_{S+a}
$$
(also denoted $\psi^R_{(m;l),a}$) represented in string diagrams as
$$
\psi^R_{S;a}=
\hbox{
\begin{tikzpicture}[baseline=(current  bounding  box.center)]
\node (m) at (0.75,0){$\underline{a}$};
\node (p) at (0,0){$P_{S}$};
\coordinate (v) at (0,0.7);
\node (pm) at (0,1.3){$P_{S+a}$};

\draw[thick,blue] (m) to [out=north,in=east] (v);
\draw[thick] (p) to (v) to (pm);
\draw[fill,blue] (v) circle [radius=0.06];
\end{tikzpicture}
},
$$
providing a basis vector for $\mathcal{P}_d(P_{S}\underline{a},P_{S+a})$ using the right unitor. Similarly to before, with analogous proof, we have:

\begin{lem}\label{rightpsiasso}
Let $a,b\in \Z_d$ and $P_S\in \mathcal{P}_d$. Then
$$
\psi^R_{S+a;b}\circ (\psi^R_{S;a} \otimes \id_{\underline{b}})= \psi^R_{S;a+b}\circ (\id_{P_{S}}\otimes\psi^R_{(a;0),b}),
$$
in string diagrams:
$$
\hbox{
\begin{tikzpicture}[baseline=(current  bounding  box.center)]

\node (a) at (0.8,0){$\underline{a}$};
\node (b) at (0.5,0){$\underline{b}$};
\node (p) at (0,0){$P_{S}$};
\coordinate (v) at (0,0.7);
\coordinate (v1) at (0,1);
\node (pm) at (0,1.6){$P_{S+a+b}$};

\draw[thick,blue] (a) to [out=north,in=east] (v1);
\draw[thick,blue] (b) to [out=north,in=east] (v);
\draw[thick] (p) to (v) to (v1) to (pm);
\draw[fill,blue] (v) circle [radius=0.06];
\draw[fill,blue] (v1) circle [radius=0.06];
\end{tikzpicture}
}
=
\hbox{
\begin{tikzpicture}[baseline=(current  bounding  box.center)]

\node (a) at (0.8,0){$\underline{a}$};
\node (b) at (0.5,0){$\underline{b}$};
\node (p) at (0,0){$P_{S}$};
\coordinate (v) at (0.5,0.6);
\coordinate (v1) at (0,1);
\node (pm) at (0,1.6){$P_{S+a+b}$};

\draw[thick,blue] (a) to [out=north,in=east] (v);
\draw[thick,blue] (b) to [out=north,in=south] (v) to [out=north, in=east] (v1);
\draw[thick] (p) to (v1) to (pm);
\draw[fill,blue] (v) circle [radius=0.06];
\draw[fill,blue] (v1) circle [radius=0.06];
\end{tikzpicture}
}.
$$
\end{lem}

These morphisms interact well with the associators in the sense that:
\begin{lem}\label{psileftrightorder}
Let $a,b\in \Z_d$ and $P_{S}\in \mathcal{P}_d$. Then
$$
\psi^R_{S+a;b}\circ (\psi^L_{a;S} \otimes \id_{\underline{b}})= \psi^L_{S+b;a}\circ (\id_{\underline{a}}\otimes\psi^R_{S;b}),
$$
which in string diagrams is
$$
\hbox{
\begin{tikzpicture}[baseline=(current  bounding  box.center)]
\node (a) at (-0.75,0){$\underline{a}$};
\node (b) at (0.75,0){$\underline{b}$};
\node (p) at (0,0){$P_{S}$};
\coordinate (v) at (0,0.6);
\coordinate (v1) at (0,0.8);
\node (pm) at (0,1.3){$P_{S+a+b}$};

\draw[thick,blue] (a) to [out=north,in=west] (v);
\draw[thick,blue] (b) to [out=north,in=east] (v1);
\draw[thick] (p) to (v) to (v1) to (pm);
\draw[fill,blue] (v) circle [radius=0.06];
\draw[fill,blue] (v1) circle [radius=0.06];
\end{tikzpicture}
}
=
\hbox{
\begin{tikzpicture}[baseline=(current  bounding  box.center)]
\node (a) at (-0.75,0){$\underline{a}$};
\node (b) at (0.75,0){$\underline{b}$};
\node (p) at (0,0){$P_{S}$};
\coordinate (v) at (0,0.8);
\coordinate (v1) at (0,0.6);
\node (pm) at (0,1.3){$P_{S+a+b}$};

\draw[thick,blue] (a) to [out=north,in=west] (v);
\draw[thick,blue] (b) to [out=north,in=east] (v1);
\draw[thick] (p) to (v) to (v1) to (pm);
\draw[fill,blue] (v) circle [radius=0.06];
\draw[fill,blue] (v1) circle [radius=0.06];
\end{tikzpicture}
}.
$$

\end{lem}
\begin{proof}
This is immediate from the naturality of the unitors and the additivity of $s$ (Equation \eqref{sisadditive}.
\end{proof}

Using Lemma \ref{sandunitorinteraction}, we have an alternative form for these maps
\begin{align*}
    \psi^L_{a;S}&= s_{\mi a,0}^{-1} \circ \,_{\mi a}\lambda \circ (s_{\mi a, 0}\otimes \id_{P_{S}})\\
    \psi^R_{S;a}&= s_{0,\mi a}^{-1} \circ \rho_{\mi a} \circ (\id_{P_{S}} \otimes s_{0,\mi a}).
\end{align*}
The underlying modules for $\underline{a}P_{S}$ (and $P_{S}\underline{a}$) are sums of $\C[x,y,z]$, that we can index by $i,j\in \{0,1\}$ corresponding to $\underline{a}_i\left(P_{m;l}\right)_j$. Recall from Equation \eqref{unitors} that the unitors are zero on the 1-summands of the unit. Tracing through the maps we find that the only nonzero components (in the sense of Notation \ref{gradingcomponentsnotation}) are
\begin{equation}
    \begin{split}
        \label{psilexplicit}    \left(\psi^L_{a;S}\right)_{0j}\colon& f(x,y,z) \mapsto \eta^{-a \vert S \vert j} f(x, \eta^{-a}x,z) \\
    \left(\psi^R_{S;b}\right)_{i0}\colon& f(x,y,z) \mapsto  f(x, \eta^{b}z,z).
    \end{split}
\end{equation}

The maps $\psi^L$ and $\psi^R$ give us a priori two different choices of basis for $\mathcal{P}_d(\underline{a}\underline{b},\underline{a+b})$, we will show they actually agree. This is in a sense a generalization to the result that $\lambda_{\mathbb{I}}=\rho_{\mathbb{I}}$ \cite{Carqueville2010}. In $\mathrm{HMF}_{\mathrm{bi}} \left( W \right)$ and consequently also in $\mathcal{P}_d$, morphisms are taken up to homology (see Definition \ref{hmfdef}), and so to prove this equality we need to find (\`a la \cite[Appendix A]{Carqueville2010}) a degree 1 morphism $h$ such that
\begin{equation}\label{homotopyunitordef}
\lambda_\mathbb{I}-\rho_{\mathbb{I}}=d^{\mathbb{I}} h+ h d^{\mathbb{II}},    
\end{equation}
with components $h_{ij}$, where $i,j\in \{0,1\}$ refers to the summand $\mathbb{I}_i \otimes\mathbb{I}_j$ in $\mathbb{II}$, a slight variation on Notation \ref{gradingcomponentsnotation} for odd morphisms.

\begin{lem}\label{psisymmetry}
For $a,b\in \Z_d$ we have, in $\mathcal{P}_d(\underline{a}\underline{b}, \underline{a+b})$,
$$
\psi^L_{a,(b;0)}= \psi^R_{(a;0),b},
$$
so
$$
\hbox{
\begin{tikzpicture}[baseline=(current  bounding  box.center)]
\node (m) at (-0.5,0){$\underline{a}$};
\node (p) at (0,0){$\underline{b}$};
\coordinate (v) at (0,0.7);
\node (pm) at (0,1.3){$\underline{a+b}$};

\draw[thick,blue] (m) to [out=north,in=west] (v);
\draw[thick,blue] (p) to (v) to (pm);
\draw[fill,blue] (v) circle [radius=0.06];
\end{tikzpicture}
}
=
\hbox{
\begin{tikzpicture}[baseline=(current  bounding  box.center)]
\node (m) at (0.5,0){$\underline{b}$};
\node (p) at (0,0){$\underline{a}$};
\coordinate (v) at (0,0.7);
\node (pm) at (0,1.3){$\underline{a+b}$};

\draw[thick,blue] (m) to [out=north,in=east] (v);
\draw[thick,blue] (p) to (v) to (pm);
\draw[fill,blue] (v) circle [radius=0.06];
\end{tikzpicture}
}.
$$
\end{lem}

\begin{proof}
We will show that the morphism $\tilde{h}$ given by $\tilde{h}_{00}=\eta^{-a}h_{00}$, $\tilde{h}_{10}=\eta^{a}h_{10}$, $\tilde{h}_{01}=h_{01}$, and $\tilde{h}_{11}=h_{11}$ for $h$ as in Equation \eqref{homotopyunitordef} defines a chain homotopy between $\psi^L_{a,(b;0)}$ and $\psi^R_{(a;0),b}$. That is, we want to show that $$
\psi^L_{a,(b;0)}- \psi^R_{(a;0),b}= d^{\underline{a+b}} \tilde{h}+ \tilde{h}d^{\underline{a}\underline{b}}.
$$

Recall that we have
\begin{align*}
    d^a_1(x,y) &= x-\eta^a y \\
    d^a_0(x,y) &= \prod_{i\neq a}(x - \eta^i y)
\end{align*}
as the differentials with domain $\underline{a}_1$ and $\underline{a}_0$, respectively. To unclutter the notation, we will write $d_i$ for $d^\mathbb{I}_i$ in the rest of this proof.

On $\underline{a}_0\underline{b}_0$ we need, for $f\in \C[x,y,z]$,
\begin{align*}
f(x,\eta^{-a}x, z) -f(x, \eta^b z, z)=& d^{a+b}_1 \eta^{-a} h_{00}(f(x,y,z)) + \eta^{a}h_{10}(d^a_0(x,y)f(x,y,z)\\ &+ h_{01}(d^b_{0}(y,z)f(x,y,z)).
\end{align*}
Without loss of generality, we can send $x \mapsto \eta^a x$, $z \mapsto \eta^{-b} z$ and replace $f$ by $f(\eta^{-a}\cdot, \cdot, \eta^{-b}\cdot)$. Note that this acts on the differentials as
\begin{align*}
    d^{a}_1(\eta^a x, y)=\eta^a d_1( x,y ) &\quad& d^{a}_0(\eta^a x, y)=\eta^{-a} d_1(x,y)\\
    d^{b}_1 (y,\eta^{-b}z) = d_{1}(y,z) &\quad& d^{b}_0 (y,\eta^{-b}z) = d_{0}(y,z)\\
    d_{1}^{a+b}(\eta^a x, \eta^{-b}z) = \eta^a d_1(x,z) &\quad& d_{0}^{a+b}(\eta^a x, \eta^{-b}z) = \eta^{-a} d_0(x,y).
\end{align*}
This means we get, using that $h$ is $\C$-linear,
\begin{align*}
f(x,x,z)-f(x,z,z)=&d_1(x,z) h_{00} f(x,y,z) + h_{10}(d_0(x,y)f(x,y,z))\\ &+ h_{01}(d_0(y,z)f(x,y,z)),
\end{align*}
and this is exactly the $00$ component of Equation \eqref{homotopyunitordef}. Using the same transformation, one similarly shows that the other components hold.
\end{proof}

The maps $\psi^L$ (or equivalently $\psi^R$) give maps
$$
\ev'_{\underline{a}} := \psi^L_{\mi a,(a,0)}\colon \underline{\mi a}\underline{a} \rar \mathbb{I},
$$
that we can use to expose $\underline{\mi a}$ as the left dual of $\underline{a}$. Suppressing as is customary unit strands from the string diagrams, this means we set
$$
\hbox{
\begin{tikzpicture}[baseline=(current  bounding  box.center)]
\node (m) at (-0.5,0){$\underline{\mi a}$};
\node (p) at (0,0){$\underline{ a}$};
\coordinate (v) at (-0.25,0.5);

\draw[thick,blue] (m) to [out=north,in=west] (v) to [out=east,in=north] (p);
\end{tikzpicture}
}
=
\hbox{
\begin{tikzpicture}[baseline=(current  bounding  box.center)]
\node (m) at (-0.5,0){$\underline{\mi a}$};
\node (p) at (0,0){$\underline{a}$};
\coordinate (v) at (0,0.7);

\draw[thick,blue] (m) to [out=north,in=west] (v);
\draw[thick,blue] (p) to (v);
\draw[fill,blue] (v) circle [radius=0.06];
\end{tikzpicture}
}.
$$

We claim that this evaluation agrees with the evaluation from Equation \ref{evforpermtype}. Because $\underline{a}$ is invertible any two choices of evaluation are related by a scalar multiplication and hence it suffices to show that $\ev'_{\underline{a}}\circ \coev_{\underline{-a}}=\id_{\mathbb{I}}$:

\begin{lem}\label{psievagree}
Let $a\in \cat{V}_d$ and take $\ev_{\underline{a}}'$ as defined above, and $\coev_{\underline{\mi a}}$ from Equation \ref{coevforpermtype}. Then
$$
\ev_{\underline{a}}' \circ \coev_{\underline{\mi a}}= \id_\mathbb{I}.
$$
\end{lem}
\begin{proof}
We work degree by degree, using Notation \ref{gradingcomponentsnotation}. On the degree 1 summand $\mathbb{I}_{1}=\C[x,z]$ we need to compute
$$
 (\ev_{a}')\circ(\coev_{\mi a})_{1}=  (\ev_{a}')_{10}\circ(\coev_{\mi a})_{10}+  (\ev_{a}')_{01}\circ(\coev_{\mi a})_{01}.
$$
Now $(\ev'_{a})_{10}\circ(\coev_{\mi a})_{10}  =0$ as $(\ev'_{a})_{10}=0$, and using that $(\coev_{\mi a})_{01}=\eta^{\mi a}$ we see
$$
 (\ev'_{a})_{01}\circ(\coev_{\mi a})_{01} (f(x,z)) =   \eta^{a}\eta^{\mi a} f(x,z)= f(x,z),
$$
where the factor $\eta^{a}$ comes from the evaluation, see Equation \ref{psilexplicit}. So the degree 1 part of this map is indeed the identity.

For the degree 0 component, we have $(\ev'_{a})_{11}=0$, leaving us with
\begin{align*}
    (\ev'_{a})\circ(\coev_{\mi a})_{0}(f(x,z))   &= (\ev_{a})_{00}\circ(\coev_{\mi a})_{00} (f(x,z))\\
                                                &= (\ev_{a})_{00}(f(x,z))\\
                                                &= f(x,z).
\end{align*}
This finishes the proof.
\end{proof}

We will denote the coevaluation in string diagrams as
$$
\coev_{a}= \hbox{
\begin{tikzpicture}[baseline=(current  bounding  box.center)]
\node (m) at (-0.5,0){$\underline{a}$};
\node (p) at (0,0){$\underline{\mi a}$};
\coordinate (v) at (-0.25,-0.5);

\draw[thick,blue] (m) to [out=south,in=west] (v) to [out=east,in=south] (p);
\end{tikzpicture}
}
$$

We will now move on to considering the interaction of these basis vectors with fusion of arbitrary simples. In order to do this, we need some information about the fusion isomorphisms from \cite[Appendix B]{Davydov2018}. These are inductively defined basis vectors for those $\mathcal{P}_d(P_{m;l},P_{m';l'}P_{m'';l''})$ that are nontrivial, starting from
\begin{equation} 
\begin{split}
g^-_{m,(m';l')} & \colon P_{m+m'+1;l'-1} \rar P_{m;1}P_{m';l'}, \quad\tn{and} \\ \quad g^+_{m,(m';l')} & \colon P_{m+m';l'+1} \rar P_{m;1}P_{m';l'}.   
\end{split}
\label{gplusminus}
\end{equation}
We will not need the explicit form of these maps. It will suffice to know that as maps $(P_{m+m'+\frac{1}{2}(1-\epsilon);l'+\epsilon})_{i+j}\cong \C[x,y]\rar (P_{m;1})_i(P_{m';l'})_j \cong \C[x,y,z]$ they are determined by multiplication by homogeneous $g_{ij}^\epsilon\in \C[x,y,z]$. Writing $\epsilon$ for $\pm1$, the degrees of these polynomials are
\begin{align}
\deg((g_{m,(m';l')}^\epsilon)_{ij})=(-1)^j il+ij(d-2)+\frac{(-1)^i}{2} \left( 1-\left( -1 \right)^j \epsilon \right)
        \label{degreesg}
\end{align}
for $i,j \in \lbrace 0,1 \rbrace$. In what follows we will depict these inductively defined morphisms $\varphi^{S',S''}_{S}\in \mathcal{P}_d(P_{S},P_{S'}P_{S''})$ 
simply by a trivalent vertex,
$$
\varphi^{S',S''}_{S}=
\hbox{
\begin{tikzpicture}[baseline=(current  bounding  box.center)]
\node (p1) at (0,-0.2){$P_{S}$};
\node (p2) at (-0.5,1){$P_{S'}$};
\node (p3) at (0.5,1){$P_{S''}$};

\coordinate (v) at (0,0.4);
\draw[fill] (v) circle [radius=0.06];

\draw[thick] (p1) to (v);
\draw[thick] (p2) to [out=south, in=west] (v) to [out=east,in=south] (p3);
\end{tikzpicture}
}
$$
for $|S'|,|S''|>1$. Note that 
$$
\varphi^{(m;1),(m';l')}_{(m+m',l'+\epsilon)}=g^\epsilon_{m,(m';l')}.
$$

The reader will notice that these $\varphi^{S,S'}_{S''}$ map into the tensor product $P_{S'}P_{S''}$ rather than out of it as $\psi^L$ and $\psi^R$ do. We will therefore work with $(\psi^L)^{-1}$ and $(\psi^R)^{-1}$ instead\footnote{The reader may wonder why we introduced $\psi^L$ rather than $(\psi^L)^{-1}$ in the first place. This is because the unitors are easier to describe than their inverses, but the basis vectors for $\mathcal{P}_d(P_{m;l},P_{m';l'}P_{m'';l''})$ can be expressed as multiplication by a polynomial while their inverses are more complicated.}, and we will depict them by
$$
(\psi^L_{a,S})^{-1}=
\hbox{
\begin{tikzpicture}[baseline=(current  bounding  box.center)]
\node (m) at (-0.75,0){$\underline{a}$};
\node (p) at (0,0){$P_{S}$};
\coordinate (v) at (0,-0.7);
\node (pm) at (0,-1.3){$P_{S+a}$};

\draw[thick,blue] (m) to [out=south,in=west] (v);
\draw[thick] (p) to (v) to (pm);
\draw[fill,blue] (v) circle [radius=0.06];
\end{tikzpicture}
}
\quad\mbox{and}\quad
(\psi^R_{S,a})^{-1}=
\hbox{
\begin{tikzpicture}[baseline=(current  bounding  box.center)]
\node (m) at (0.75,0){$\underline{a}$};
\node (p) at (0,0){$P_{S}$};
\coordinate (v) at (0,-0.7);
\node (pm) at (0,-1.3){$P_{S+a}$};

\draw[thick,blue] (m) to [out=south,in=east] (v);
\draw[thick] (p) to (v) to (pm);
\draw[fill,blue] (v) circle [radius=0.06];
\end{tikzpicture}
}.
$$
Observe that, because they only involve isomorphisms, Lemmas \ref{leftpsiasso}, \ref{rightpsiasso}, \ref{psileftrightorder}, and \ref{psisymmetry} have clear analogues for these inverses that are obtained by reading the string diagrams downwards.

The result we are after is:

\begin{prop}\label{conjugationtrivfusion}
Let $a\in \cat{V}_d$ and $P_{m;l}, P_{m';l'},P_{m'';l''}\in \mathcal{P}_d$. Then
$$
\hbox{
\begin{tikzpicture}[baseline=(current  bounding  box.center)]
\node (a) at (-1.4,2){$\underline{a}$};
\node (am) at (1.4,2){$\underline{\mi a}$};
\node (pl) at (-0.5,2){$P_{m';l'}$};
\node (pr) at (0.5,2){$P_{m'';l''}$};
\coordinate (v1) at (0,1.2);
\coordinate (v2) at (0,0.8);
\coordinate (v3) at (0,0.6);
\node (pm) at (0,0){$P_{m;l}$};

\draw[thick,blue] (a) to [out=south,in=west] (v2);
\draw[thick,blue] (am) to [out=south,in=east] (v3);
\draw[thick] (pl) to  [out=south, in=west] (v1) to [out=east, in=south] (pr);
\draw[thick] (pm) to (v1);
\draw[fill] (v1) circle [radius=0.05];
\draw[fill,blue] (v2) circle [radius=0.05];
\draw[fill,blue] (v3) circle [radius=0.05];
\end{tikzpicture}
}
=
\eta^{\frac{a}{2}(l-l'-l'')}
\hbox{
\begin{tikzpicture}[baseline=(current  bounding  box.center)]
\node (a) at (-1.4,2){$\underline{a}$};
\node (am) at (1.4,2){$\underline{\mi a}$};
\node (pl) at (-0.5,2){$P_{m';l'}$};
\node (pr) at (0.5,2){$P_{m'';l''}$};
\coordinate (v1) at (0,0.7);
\coordinate (v2l) at (-0.5,1);
\coordinate (v2r) at (0.5,1);
\coordinate (v3l) at (-0.5,1.2);
\coordinate (v3r) at (0.5,1.2);
\coordinate (c) at (0,1.4);
\node (pm) at (0,0){$P_{m;l}$};

\draw[thick,blue] (a) to [out=south,in=west] (v3l);
\draw[thick,blue] (am) to [out=south,in=east] (v3r);
\draw[thick] (pl) to (v3l) to [out= south, in=north] (v2l) to [out=south, in=west] (v1) to [out=east, in=south] (v2r) to  (v3r) to (pr);
\draw[thick] (pm) to (v1);
\draw[thick,blue] (v2l) to [out =east, in=west] (c) to [out=east, in=west] (v2r);
\draw[fill] (v1) circle [radius=0.06];
\draw[fill,blue] (v2l) circle [radius=0.06];
\draw[fill,blue] (v2r) circle [radius=0.06];
\draw[fill,blue] (v3r) circle [radius=0.06];
\draw[fill,blue] (v3l) circle [radius=0.06];
\end{tikzpicture}
},
$$
where we have represented $\ev_{a}\colon \underline{-a}\underline{a}\rar \mathbb{I}$ by a cap.
\end{prop}

Whenever $P_{m;l}$ does not appear as a summand of $P_{m';l'}P_{m'';l''}$, this proposition is trivially true. We will proceed by induction on $l'$, mirroring the argument from \cite[Appendix B]{Davydov2018}. Note that picking $m',l',m'',l''$ determines the possible choices for $m,l$. The base case for the induction ($l'=1$) is:

\begin{lem}
For $a\in \cat{V}_d$ and the $g^\epsilon_{m',(m;l)}$ from Equation \eqref{gplusminus}, we have
$$
\hbox{
\begin{tikzpicture}[baseline=(current  bounding  box.center)]
\node (a) at (-1.4,2){$\underline{a}$};
\node (am) at (1.4,2){$\underline{\mi a}$};
\node (pl) at (-0.5,2){$P_{m';1}$};
\node (pr) at (0.5,2){$P_{m'';l''}$};
\coordinate (v1) at (0,1.2);
\coordinate (v2) at (0,0.8);
\coordinate (v3) at (0,0.6);
\node (pm) at (0,0){$P_{m'+m''+\frac{1}{2}(1-\epsilon);l''+ \epsilon}$};

\draw[thick,blue] (a) to [out=south,in=west] (v2);
\draw[thick,blue] (am) to [out=south,in=east] (v3);
\draw[thick] (pl) to  [out=south, in=west] (v1) to [out=east, in=south] (pr);
\draw[thick] (pm) to (v1);
\draw[fill] (v1) circle [radius=0.05];
\draw[fill,blue] (v2) circle [radius=0.05];
\draw[fill,blue] (v3) circle [radius=0.05];
\end{tikzpicture}
}
=
\eta^{\frac{a}{2}(\epsilon-1)}
\hbox{
\begin{tikzpicture}[baseline=(current  bounding  box.center)]
\node (a) at (-1.4,2){$\underline{a}$};
\node (am) at (1.4,2){$\underline{\mi a}$};
\node (pl) at (-0.5,2){$P_{m';1}$};
\node (pr) at (0.5,2){$P_{m'';l''}$};
\coordinate (v1) at (0,0.7);
\coordinate (v2l) at (-0.5,1);
\coordinate (v2r) at (0.5,1);
\coordinate (v3l) at (-0.5,1.2);
\coordinate (v3r) at (0.5,1.2);
\coordinate (c) at (0,1.4);
\node (pm) at (0,0){$P_{m'+m''+\frac{1}{2}(1-\epsilon);l''+ \epsilon}$};

\draw[thick,blue] (a) to [out=south,in=west] (v3l);
\draw[thick,blue] (am) to [out=south,in=east] (v3r);
\draw[thick] (pl) to (v3l) to [out= south, in=north] (v2l) to [out=south, in=west] (v1) to [out=east, in=south] (v2r) to  (v3r) to (pr);
\draw[thick] (pm) to (v1);
\draw[thick,blue] (v2l) to [out =east, in=west] (c) to [out=east, in=west] (v2r);
\draw[fill] (v1) circle [radius=0.06];
\draw[fill,blue] (v2l) circle [radius=0.06];
\draw[fill,blue] (v2r) circle [radius=0.06];
\draw[fill,blue] (v3r) circle [radius=0.06];
\draw[fill,blue] (v3l) circle [radius=0.06];
\end{tikzpicture}
}.
$$
\end{lem}
\begin{proof}
For brevity, write $l_\epsilon= l'' +\epsilon$ and $m_\epsilon=m'+m''+\frac{1}{2}(1-\epsilon)$. We will proceed in two steps. To begin, we will show 
$$
\hbox{
\begin{tikzpicture}[baseline=(current  bounding  box.center)]
\node (a) at (-1.4,2){$\underline{a}$};
\node (am) at (1.4,2){$\underline{\mi a}$};
\node (pl) at (-0.5,2){$P_{m';1}$};
\node (pr) at (0.5,2){$P_{m'';l''}$};
\coordinate (v1) at (0,1.2);
\coordinate (v2) at (0,0.8);
\coordinate (v3) at (0,0.6);
\node (pm) at (0,0){$P_{m_\epsilon;l_\epsilon}$};

\draw[thick,blue] (a) to [out=south,in=west] (v2);
\draw[thick,blue] (am) to [out=south,in=east] (v3);
\draw[thick] (pl) to  [out=south, in=west] (v1) to [out=east, in=south] (pr);
\draw[thick] (pm) to (v1);
\draw[fill] (v1) circle [radius=0.05];
\draw[fill,blue] (v2) circle [radius=0.05];
\draw[fill,blue] (v3) circle [radius=0.05];
\end{tikzpicture}
}
=
\eta^{\frac{a}{2}(\epsilon-1)}
\hbox{
\begin{tikzpicture}[baseline=(current  bounding  box.center)]
\node (a) at (-1.4,2){$\underline{a}$};
\node (am) at (1.4,2){$\underline{\mi a}$};
\node (pl) at (-0.5,2){$P_{m';1}$};
\node (pr) at (0.5,2){$P_{m'';l''}$};
\coordinate (v1) at (0,0.7);
\coordinate (v2l) at (-0.5,1);
\coordinate (v2r) at (0.5,1);
\coordinate (v3l) at (-0.5,1.2);
\coordinate (v3r) at (0.5,1.2);
\coordinate (c) at (0,1.4);
\node (pm) at (0,0){$P_{m_\epsilon;l_\epsilon}$};

\draw[thick,blue] (a) to [out=south,in=west] (v3l);
\draw[thick,blue] (am) to [out=south,in=east] (v3r);
\draw[thick] (pl) to (v3l) to [out= south, in=north] (v2l) to [out=south, in=west] (v1) to [out=east, in=south] (v2r) to  (v3r) to (pr);
\draw[thick] (pm) to (v1);
\draw[fill] (v1) circle [radius=0.06];
\draw[fill,blue] (v2l) circle [radius=0.06];
\draw[fill,blue] (v2r) circle [radius=0.06];
\draw[fill,blue] (v3r) circle [radius=0.06];
\draw[fill,blue] (v3l) circle [radius=0.06];
\end{tikzpicture}
},
$$
where the left and right non-vertex dots in the right hand diagram represent $s_{0,-a}$ and $s_{a,0}$, respectively. The non-trivial step is showing that the upper and lower composites in
\begin{center}
\begin{tikzcd}
& \,_a(P_{m_\epsilon;l_\epsilon})_{\mi a} \arrow[dr,"\,_a(g^\epsilon_{m'',(m';l')})_{\mi a}"] & \\
P_{m'';l}  \arrow[ur,"s_{a,-a}"] \arrow[dr,"g^\epsilon_{m'',(m';l')}"'] & & \,_a(P_{m';1})_{\mi a} \otimes_a P_{m'';l''})_{\mi a}\\
& _{a}\left( P_{m';1} \otimes P_{m'';l''} \right)_{-a} \arrow[ur,"s_{a,-a} \cdot s_{a,-a}"']&
\end{tikzcd},
\end{center}
coming from the left and right hand diagrams, respectively, agree up to $\eta^{\frac{a}{2}(1-\epsilon)}$. Because the polynomials $(g^\epsilon_{m'',(m';l')})_{ij}$ are homogeneous, we have that $\sigma_x^{-a} \otimes \sigma_y^{-a} \otimes \sigma_z^{-a} \left( g^\epsilon_{m'',(m';l')} \right)=\eta^{-a\mathrm{deg} (g^\epsilon_{m'',(m';l')})} g^\epsilon_{m',(m;l)}$, where $\mathrm{deg} \left( g^\epsilon_{m'',(m';l')} \right)$ denotes the degree of the maps $g^\epsilon_{m'',(m';l')}$, see Equation \eqref{degreesg}. Now, consider an element $f\left(x,z \right)$ in $ P_{m_\epsilon;l_\epsilon} $ of degree $k$, where $k=0$ if $i=j$ and $k=1$ otherwise. One can rewrite this as $k=\frac{1}{2}\left(1-\left(-1\right)^{i+j}\right)$. Tracing this element through the composition of maps we obtain 
$$
g_{ij}^\epsilon(x,y,z) \eta^{-\frac{ a}{2}(l+\epsilon+1)(1-(-1)^{i+j})} f \left(\eta^{\mi a}x,\eta^{\mi a}z \right)
$$
for the upper composite, and
$$
\eta^{-a(2i+(l+1)j+ \mathrm{deg}(g^\epsilon_{m'',(m';l')})} g_{ij}^\epsilon(x,y,z) f \left(\eta^{\mi a} x,\eta^{\mi a} z \right)
$$
for the bottom route. These indeed differ by a factor of $\eta^{\frac{a}{2} (\epsilon-1)}$.

What remains to be checked is the equality:
$$
\hbox{
\begin{tikzpicture}[baseline=(current  bounding  box.center)]
\node (pl) at (-0.6,2){$P_{m'+a;1}$};
\node (pr) at (0.6,2){$P_{m''\mi a;l''}$};
\coordinate (v1) at (0,1.2);
\coordinate (vl) at (-0.6,1.5);
\coordinate (vr) at (0.6,1.5);
\node (pm) at (0,0){$P_{m_\epsilon;l_\epsilon}$};

\draw[thick] (pl) to [out=south,in=north] (vl) to [out=south, in=west] (v1) to [out=east, in=south] (vr) to [out=north,in=south] (pr);
\draw[thick] (pm) to (v1);
\draw[fill] (v1) circle [radius=0.05];
\draw[fill,blue] (vl) circle [radius=0.05];
\draw[fill,blue] (vr) circle [radius=0.05];
\end{tikzpicture}
}
=
\hbox{
\begin{tikzpicture}[baseline=(current  bounding  box.center)]

\node (pl) at (-0.6,2){$P_{m'+a;1}$};
\node (pr) at (0.6,2){$P_{m''\mi a;l''}$};
\coordinate (v1) at (0,0.7);
\coordinate (v2l) at (-0.6,1);
\coordinate (v2r) at (0.6,1);
\coordinate (v3l) at (-0.6,1.2);
\coordinate (v3r) at (0.6,1.2);
\coordinate (c) at (0,1.4);
\node (pm) at (0,0){$P_{m_\epsilon;l_\epsilon}$};

\draw[thick] (pl) to (v3l) to [out= south, in=north] (v2l) to [out=south, in=west] (v1) to [out=east, in=south] (v2r) to  (v3r) to (pr);
\draw[thick] (pm) to (v1);
\draw[thick,blue] (v2l) to [out =east, in=west] (c) to [out=east, in=west] (v2r);
\draw[fill] (v1) circle [radius=0.06];
\draw[fill,blue] (v2l) circle [radius=0.06];
\draw[fill,blue] (v2r) circle [radius=0.06];

\end{tikzpicture}
}
$$
The inverses of $\psi^{L}$ and $\psi^R$ on the right hand side start with $s_{0,-a}$ and $s_{a,0}$, just like appear in the left hand diagram. The next maps are then inverse unitors $\rho_{(P_{m';1})_{a}}^{-1}\otimes \lambda^{-1}_{\,_{\mi a} P_{m'',l''}}$. By the triangle equations for the unitors, we can replace this by $\lambda_{\mathbb{I}}^{-1}\otimes \id_{\,_{\mi a}P_{m'';l''}} \circ \lambda^{-1}_{\,_{\mi a} P_{m'',l''}}$. This is then followed by
\begin{equation*}
\,_a\mathbb{I} \mathbb{I}_{\mi a}\xrightarrow{ s_{a,0}^{-1}\cdot s_{0,\mi a}^{-1}} \underline{\mi a} \underline{a} \xrightarrow{\ev_{a}} \mathbb{I}.
\end{equation*}
Using Lemma \ref{psievagree}, we see that this is the same as
\begin{align*}
	\lambda_{\mathbb{I}}\circ (\id_{\mathbb{I}}\otimes (s_{a,0})^{-1})\circ (s_{0,a}\otimes \id_{\underline{a}})\circ(s_{a,0}^{-1}\otimes s_{0,\mi a}^{-1}&)=\lambda_{\mathbb{I}}\circ(s_{0,a}s_{a,0}^{-1}\otimes s^{-1}_{a,a})\\
					&= s_{\mi a,a}^{-1}\circ \lambda_{\,_{\mi a}\mathbb{I}_a}\circ(s_{\mi a, a})\\
					&= s_{\mi a,a}^{-1}\circ \,_{\mi a}\lambda_{\mathbb{I}_a},
\end{align*}
where the last equality is Lemma \ref{sandunitorinteraction}. Composing this with the inverse unitors we get
$$
 \,_{\mi a}\lambda_{\mathbb{I}_a}\circ (\id_{\mathbb{I}}\otimes s_{\mi a,a}^{-1}\otimes \id_{\,_{\mi a}P_{m'';l''}} )\circ (\lambda_{\mathbb{I}}^{-1}\otimes \id_{\,_{\mi a}P_{m'';l''}}) \circ \lambda^{-1}_{\,_{\mi a} P_{m'',l''}},
$$
where we have omitted some identity morphisms for brevity. Using Lemma \ref{sandunitorinteraction} on $s_{\mi a,a}^{-1}$ and $\lambda_{\mathbb{I}}^{-1}$, we see this is just a composition of unitors, and we are done.
\end{proof}

The induction step is then:

\begin{proof}[Proof of Proposition \ref{conjugationtrivfusion}]
Assume that Proposition \ref{conjugationtrivfusion} holds for all $a,m,m',m''\in\Z_d$, $l,l''\in \{0,1, \dots, d-2\}$ and $l'<L$, we want to show that it holds for $l'=L+1$. For our induction we will use that $P_{0;1}P_{m;L}$ includes along $g^{+}_{0,(m;L)}$ onto $P_{m;L+1}$, and this allows us to use that the proposition holds for $P_{0;1}$ and $P_{m;L}$. Using the associators we have for any $m,m',m''$ and $l',l''$ that
$$
(\psi^{(0;1),(m;L)}_{(m;L+1)}\otimes \id_{P_{m';l'}})\psi^{(m;L+1),(m';l')}_{(m'',l'')}=\alpha (\id_{P_{0;1}}\otimes \psi^{(m;L),(m';l')}_{(m';l''-1)})\psi_{(m'';l'')}^{(0;1),(m';l''-1)}
$$ 
for some $\alpha\in \C^*$. This allows us to compute
\begin{eqnarray*}
& \hbox{
\begin{tikzpicture}[baseline=(current  bounding  box.center)]
\node (a) at (-2,2.5){$\underline{a}$};
\node (am) at (2,2.5){$\underline{\mi a}$};
\node (pl) at (-1.2,2.5){$P_{0;1}$};
\node (p) at (0,2.5){$P_{m;L}$};
\node (pr) at (1.2,2.5){$P_{m';l'}$};
\coordinate (v1) at (0,1.2);
\coordinate (v2) at (0,0.8);
\coordinate (v3) at (0,0.6);
\coordinate (v4) at (-0.6,1.9);
\node (pm) at (0,0){$P_{m'';l''}$};

\draw[thick,blue] (a) to [out=south,in=west] (v2);
\draw[thick,blue] (am) to [out=south,in=east] (v3);
\draw[thick] (v4) to  [out=south, in=west] (v1) to [out=east, in=south] (pr);
\draw[thick] (pl) to [out=south,in=west] (v4) to [out=east, in=south] (p);
\draw[thick] (pm) to (v1);
\draw[fill] (v1) circle [radius=0.05];
\draw[fill] (v4) circle [radius=0.05];
\draw[fill,blue] (v2) circle [radius=0.05];
\draw[fill,blue] (v3) circle [radius=0.05];
\end{tikzpicture}
}
= \alpha
\hbox{
\begin{tikzpicture}[baseline=(current  bounding  box.center)]
\node (a) at (-2,2.5){$\underline{a}$};
\node (am) at (2,2.5){$\underline{\mi a}$};
\node (pl) at (-1.2,2.5){$P_{0;1}$};
\node (p) at (0,2.5){$P_{m;L}$};
\node (pr) at (1.2,2.5){$P_{m';l'}$};
\coordinate (v1) at (0,1.2);
\coordinate (v2) at (0,0.8);
\coordinate (v3) at (0,0.6);
\coordinate (v4) at (0.6,1.9);
\node (pm) at (0,0){$P_{m'';l''}$};

\draw[thick,blue] (a) to [out=south,in=west] (v2);
\draw[thick,blue] (am) to [out=south,in=east] (v3);
\draw[thick] (pl) to  [out=south, in=west] (v1) to [out=east, in=south] (v4);
\draw[thick] (p) to [out=south,in=west] (v4) to [out=east, in=south] (pr);
\draw[thick] (pm) to (v1);
\draw[fill] (v1) circle [radius=0.05];
\draw[fill] (v4) circle [radius=0.05];
\draw[fill,blue] (v2) circle [radius=0.05];
\draw[fill,blue] (v3) circle [radius=0.05];
\end{tikzpicture}
}\\
&=\alpha 
\hbox{
\begin{tikzpicture}[baseline=(current  bounding  box.center)]
\node (a) at (-1.7,2.5){$\underline{a}$};
\node (am) at (1.7,2.5){$\underline{\mi a}$};
\node (pl) at (-1,2.5){$P_{0;1}$};
\node (p) at (0,2.5){$P_{m;L}$};
\node (pr) at (1,2.5){$P_{m';l'}$};
\coordinate (v1) at (0,1.2);
\coordinate (v2) at (-1,1.6);
\coordinate (v3) at (0.6,1.6);
\coordinate (v4) at (0.6,1.9);
\coordinate (v2e) at (-1,1.7);
\coordinate (v3e) at (0.6,1.7);
\coordinate (c) at (-0.3,2);
\node (pm) at (0,0){$P_{m'';l''}$};

\draw[thick,blue] (a) to [out=south,in=west] (v2);
\draw[thick,blue] (am) to [out=south,in=east] (v3);
\draw[thick,blue] (v2e) to [out= east, in=west] (c)to [out= east, in=west] (v3e);
\draw[thick] (pl) to [out=south,in=north] (v2) to  [out=south, in=west] (v1) to [out=east, in=south] (v3) to [out=north, in = south] (v4);
\draw[thick] (p) to [out=south,in=west] (v4) to [out=east, in=south] (pr);
\draw[thick] (pm) to (v1);
\draw[fill] (v1) circle [radius=0.05];
\draw[fill] (v4) circle [radius=0.05];
\draw[fill,blue] (v2) circle [radius=0.05];
\draw[fill,blue] (v3) circle [radius=0.05];
\draw[fill,blue] (v2e) circle [radius=0.05];
\draw[fill,blue] (v3e) circle [radius=0.05];
\end{tikzpicture}
}
= \alpha\eta^{\frac{a}{2}(\mi 1+l''\mi L \mi l')}
\hbox{
\begin{tikzpicture}[baseline=(current  bounding  box.center)]
\node (a) at (-1.7,2.5){$\underline{a}$};
\node (am) at (1.7,2.5){$\underline{\mi a}$};
\node (pl) at (-1,2.5){$P_{0;1}$};
\node (p) at (0,2.5){$P_{m;L}$};
\node (pr) at (1,2.5){$P_{m';l'}$};
\coordinate (v1) at (0,1.2);
\coordinate (v2) at (-1,1.6);
\coordinate (v3) at (1,1.8);
\coordinate (v4) at (0.6,1.5);
\coordinate (v2e) at (-1,1.9);
\coordinate (v3e) at (0,1.9);
\coordinate (v3e1) at (0,2);
\coordinate (v3f) at (1,1.9);
\coordinate (c) at (-0.6,2.3);
\coordinate (c1) at (0.6,2.3);
\node (pm) at (0,0){$P_{m'';l''}$};

\draw[thick,blue] (a) to [out=south,in=west] (v2);
\draw[thick,blue] (am) to [out=south,in=east] (v3);
\draw[thick,blue] (v2e) to [out= east, in=west] (c)to [out= east, in=west] (v3e);
\draw[thick,blue] (v3e1) to [out= east, in=west] (c1)to [out= east, in=west] (v3f);
\draw[thick] (pl) to [out=south,in=north] (v2) to  [out=south, in=west] (v1) to [out=east, in = south] (v4);
\draw[thick] (p) to [out=south, in=north] (v3e) to [out=south,in=west] (v4) to [out=east, in=south] (v3) to [out=north, in=south] (pr);
\draw[thick] (pm) to (v1);
\draw[fill] (v1) circle [radius=0.05];
\draw[fill] (v4) circle [radius=0.05];
\draw[fill,blue] (v2) circle [radius=0.05];
\draw[fill,blue] (v3) circle [radius=0.05];
\draw[fill,blue] (v2e) circle [radius=0.05];
\draw[fill,blue] (v3e) circle [radius=0.05];
\draw[fill,blue] (v3e1) circle [radius=0.05];
\draw[fill,blue] (v3f) circle [radius=0.05];
\end{tikzpicture}
}\\
&=\eta^{\frac{a}{2}(l\mi L \mi 1 \mi l')}
\hbox{
\begin{tikzpicture}[baseline=(current  bounding  box.center)]
\node (a) at (-1.7,2.5){$\underline{a}$};
\node (am) at (1.7,2.5){$\underline{\mi a}$};
\node (pl) at (-1,2.5){$P_{0;1}$};
\node (p) at (0,2.5){$P_{m;L}$};
\node (pr) at (1,2.5){$P_{m';l'}$};
\coordinate (v1) at (0,1.2);
\coordinate (v2) at (-1,1.8);
\coordinate (v3) at (1,1.8);
\coordinate (v4) at (-0.6,1.5);
\coordinate (v2e) at (-1,1.9);
\coordinate (v3e) at (0,1.9);
\coordinate (v3e1) at (0,2);
\coordinate (v3f) at (1,1.9);
\coordinate (c) at (-0.6,2.3);
\coordinate (c1) at (0.6,2.3);
\node (pm) at (0,0){$P_{m'';l''}$};

\draw[thick,blue] (a) to [out=south,in=west] (v2);
\draw[thick,blue] (am) to [out=south,in=east] (v3);
\draw[thick,blue] (v2e) to [out= east, in=west] (c)to [out= east, in=west] (v3e);
\draw[thick,blue] (v3e1) to [out= east, in=west] (c1)to [out= east, in=west] (v3f);
\draw[thick] (pl) to [out=south,in=north] (v2) to  [out=south, in=west] (v4) to [out=east,in=south] (v3e) to [out=north,in=south] (p);
\draw[thick] (v4) to [out=south, in= west] (v1) to [out=east, in=south] (v3) to [out=north, in=south] (pr);
\draw[thick] (pm) to (v1);
\draw[fill] (v1) circle [radius=0.05];
\draw[fill] (v4) circle [radius=0.05];
\draw[fill,blue] (v2) circle [radius=0.05];
\draw[fill,blue] (v3) circle [radius=0.05];
\draw[fill,blue] (v2e) circle [radius=0.05];
\draw[fill,blue] (v3e) circle [radius=0.05];
\draw[fill,blue] (v3e1) circle [radius=0.05];
\draw[fill,blue] (v3f) circle [radius=0.05];
\end{tikzpicture}
}
= \eta^{\frac{a}{2}(l\mi L \mi 1 \mi l')}
\hbox{
\begin{tikzpicture}[baseline=(current  bounding  box.center)]
\node (a) at (-1.7,2.5){$\underline{a}$};
\node (am) at (1.7,2.5){$\underline{\mi a}$};
\node (pl) at (-1,2.5){$P_{0;1}$};
\node (p) at (0,2.5){$P_{m;L}$};
\node (pr) at (1,2.5){$P_{m';l'}$};
\coordinate (v1) at (0,1.2);
\coordinate (v2) at (-0.6,1.6);
\coordinate (v3) at (1,2);
\coordinate (v4) at (-0.6,1.9);
\coordinate (v2e) at (-0.6,1.7);
\coordinate (v3e) at (1,1.7);
\coordinate (c) at (0.3,2);
\node (pm) at (0,0){$P_{m'';l''}$};

\draw[thick,blue] (a) to [out=south,in=west] (v2);
\draw[thick,blue] (am) to [out=south,in=east] (v3);
\draw[thick,blue] (v2e) to [out= east, in=west] (c)to [out= east, in=west] (v3e);
\draw[thick] (pl) to [out=south,in=west] (v4) to [out=east,in=south] (p);
\draw[thick] (v4) to [out=south,in=north] (v2) to  [out=south, in=west] (v1) to [out=east, in=south] (v3e) to [out=north, in = south] (v3) to [out=north,in=south] (pr);
\draw[thick] (pm) to (v1);
\draw[fill] (v1) circle [radius=0.05];
\draw[fill] (v4) circle [radius=0.05];
\draw[fill,blue] (v2) circle [radius=0.05];
\draw[fill,blue] (v3) circle [radius=0.05];
\draw[fill,blue] (v2e) circle [radius=0.05];
\draw[fill,blue] (v3e) circle [radius=0.05];
\end{tikzpicture}
}.
\end{eqnarray*}

The first equality corresponds to applying this associator, the second and third use the induction hypothesis, the fourth is the inverse of the associator, and the last equality again uses the induction hypothesis. Composing with the projection morphism $\pi \colon P_{0,1}P_{m,L}\rar P_{m;L+1}$ associated to the inclusion $g^{+}_{0,(m;L)}$ and using that $\pi g^{+}_{0,(m;L)}=\id_{P_{m,L+1}}$ now yields the result.
\end{proof}

\subsection{$\mathcal{P}_d$ as a $\mathcal{V}_d$-module tensor category}\label{pdasmodcatsection}
We will now use the results from the previous sections to equip $\mathcal{P}_d$ with the structure of a $\mathcal{V}_d$-module tensor category. That is, we will provide what is called a central lift $\mathcal{V}_d\rar \mathcal{Z}(\mathcal{P}_d)$ for the embedding $\mathcal{V}_d\hookrightarrow \mathcal{P}_d$ from Equation \eqref{vdintopd}. Rather than doing this directly, it is convenient to instead trivialise the adjoint action of $\mathcal{V}_d$ on $\mathcal{P}_d$, as we explain in the next section.

\subsubsection{Central functors and trivialisations of adjoint actions}
In this section we recall some basics about central lifts, and record the observation that under suitable conditions central lifts correspond to trivialisations of adjoint actions. This allows us to use the results from the previous sections to give a central lift for the inclusion $\cat{V}_d\hookrightarrow \cat{P}_d$ through specifying a trivialisation of the adjoint action of $\cat{V}_d$.

\begin{df} Let $\cat{B}$ be a braided monoidal category and $\cat{C}$ monoidal category.
\begin{itemize}
    \item[-] A \textit{central functor} (in the sense of \cite{Henriques2016}) from a $\cat{B}$ to $\cat{C}$ is a braided monoidal functor $\cat{B}\rar \dcentcat{C}$ into the Drinfeld center of $\cat{C}$. 
    \item[-] A \textit{central lift} of a functor $\cat{B}\rar \cat{C}$ to a central functor is a factorisation
$\cat{B\rar \dcentcat{C}} \rar \cat{C}$ of that functor through the forgetful functor $\dcentcat{C}\rar \cat{C}$. 
\end{itemize}
\end{df}
Concretely, a lift to a central functor means that we have to pick, for every $b\in \cat{B}$, a natural isomorphism $\beta_b$ (the \textit{half-braiding}) with components
$$
\beta_{b,c}\colon bc \xrightarrow{\cong} cb, \quad\tn{depicted by }
\hbox{
\begin{tikzpicture}[baseline=(current  bounding  box.center)]
\node (bi) at (0,0){$b$};
\node (ci) at (0.5,0){$c$};
\coordinate (co) at (0,1);
\coordinate (bo) at (0.5,1);

\begin{knot}[clip width=4pt]
\strand [blue,thick] (bi) to[out=north,in=south] (bo);
\strand[thick] (ci) to [out=north,in=south] (co);
\flipcrossings{1}
\end{knot}
\end{tikzpicture}
}
$$
natural in $b$, and satisfying for all $b,b'\in\cat{B}$ and $c,c'\in\cat{C}$ the hexagon equations
\begin{align}
\beta_{b,cc'}& =(\id_c \otimes \beta_{b,c'})(\beta_{b,c}\otimes \id_c')  \quad \tn{and}  \label{lifthalfbraidcondc}\\
\beta_{bb',c}&= (\beta_{b,c}\otimes \id_{b'})(\id_{b}\otimes \beta_{b',c}), \label{lifthalfbraidcondb}
\end{align}
where we have suppressed the associators. The lift is a functor between braided monoidal categories, and we require this functor to be braided in the usual sense.

Even if $\mathcal{B}$ does not carry a braiding, it still makes sense to talk about lifts of functors $\mathcal{B}\rar \mathcal{C}$ to functors $\mathcal{B}\rar\dcentcat{C}$. In what follows, this will be useful as an intermediate step, and we will refer to this as an \textit{unbraided central lift} (or \textit{unbraided central functor}).

\begin{df}
Suppose that $\cat{B}$ is a pointed fusion category, then the \emph{adjoint action associated to the inclusion $\iota \colon\cat{B}\subset \cat{C}$} is the monoidal functor
\begin{align*}
\Ad_{\iota}\colon \pic(\cat{B})&\rar \Aut_{\otimes}\cat{C}    \\
b &\mapsto (\Ad_{b}\colon c\mapsto  b  c b^*)
\end{align*}
from the Picard groupoid of invertible objects $\pic(\cat{B})$ of $\cat{B}$ to the tensor automorphisms of $\cat{C}$. Clearly $\Ad_b \circ \Ad_{b'}\cong \Ad_{bb'}$. The monoidal coherence isomorphism $\mu$ for $\Ad_b$ is
$$
\mu_{c,c'}\colon bcb^*bbc'b^*\xrightarrow{\id_{bc} \cdot \ev_{b} \cdot \id_{cb^*}} bcc'b^*.
$$
\end{df}

We can rephrase the data of an unbraided central lifts in terms of trivialisations of the adjoint action:

\begin{lem}\label{centralliftvstriv}
Let $\iota\colon\cat{B}\hookrightarrow \cat{C}$ be an embedding of a spherical pointed fusion category into a fusion category. Then a lift of this embedding to an unbraided central functor is equivalently an isomorphism
$$
\tau\colon \Ad_\iota \cong \Id_\cat{C}
$$
in the groupoid of tensor automorphisms of $\cat{C}$.
\end{lem}

\begin{proof}
Given a central lift $\cat{B}\rar \dcentcat{C}$, define the components of $\tau$ to be the natural isomorphisms $\tau_b$ with components
\begin{equation*}
\tau_{b,c}:= b cb^* \xrightarrow{\beta_{b,c} \cdot \id_{b*}} cbb^*\xrightarrow{\id_c \cdot \ev_b} c,
\end{equation*}
or in string diagrams
$$
\hbox{
\begin{tikzpicture}[baseline=(current  bounding  box.center)]
\node (bi) at (0,0){$b$};
\node (ci) at (0.5,0){$c$};
\coordinate (co) at (0.5,1.4);
\coordinate (w) at (0,0.4);
\coordinate (e) at (1,0.4);
\node (bo) at (1,0){$b^*$};

\begin{knot}[clip width=4pt]
\strand [blue,thick] (bi) to[out=north,in=south] (w) to [out=north, in=north] (e) to (bo);
\strand[thick] (ci) to [out=north,in=south] (co);
\flipcrossings{1}
\end{knot}
\end{tikzpicture}
}
:=
\hbox{
\begin{tikzpicture}[baseline=(current  bounding  box.center)]
\node (bi) at (0,0){$b$};
\node (ci) at (0.5,0){$c$};
\coordinate (co) at (0.5,1.4);
\coordinate (w) at (0,0.3);
\coordinate (e) at (1,0.6);
\node (bo) at (1,0){$b^*$};

\begin{knot}[clip width=4pt]
\strand [blue,thick] (bi) to[out=north,in=south] (w) to [out=north, in=north] (e) to (bo);
\strand[thick] (ci) to [out=north,in=south] (co);
\flipcrossings{1}
\end{knot}
\end{tikzpicture}
}
.
$$
This is natural in $c$ by naturality of $\beta_b$. These components are indeed monoidal in $c$, use Equation \eqref{lifthalfbraidcondc} to see that the two compositions in the diagram 
\begin{center}
\begin{tikzcd}
bcb^*bc'b^* \arrow[rd,"\tau_{b,c} \cdot \tau_{b,c'}"]& \\
bcc'b^*\arrow[u,"\mu_{c,c'}^{-1}"]\arrow[r,"\tau_{b,cc'}"']&cc'
\end{tikzcd}
\end{center}
are
$$
\hbox{
\begin{tikzpicture}[baseline=(current  bounding  box.center)]
\node (bi) at (0,0){$b$};
\node (ci) at (0.5,0){$c$};
\coordinate (co) at (0.5,1.4);
\coordinate (w) at (0,0.3);
\coordinate (e) at (1,0.6);
\coordinate (w1) at (1.3,0.3);
\node (ci1) at (1.8,0){$c'$};
\coordinate (co1) at (1.8,1.4);
\coordinate (e1) at (2.3,0.6);
\node (bo1) at (2.3,0){$b^*$};

\begin{knot}[clip width=4pt]
\strand [blue,thick] (bi) to[out=north,in=south] (w) to [out=north, in=north] (e) to [out=south,in=south] (w1) to [out=north,in=north] (e1) to (bo1);
\strand[thick] (ci) to [out=north,in=south] (co);
\strand[thick] (ci1) to [out=north,in=south] (co1);
\flipcrossings{1,2}
\end{knot}
\end{tikzpicture}
}
=
\hbox{
\begin{tikzpicture}[baseline=(current  bounding  box.center)]
\node (bi) at (0,0){$b$};
\node (ci) at (0.5,0){$c$};
\coordinate (co) at (0.5,1.4);
\coordinate (w) at (0,0.3);
\coordinate (e) at (1.3,0.6);
\node (ci1) at (0.8,0){$c'$};
\coordinate (co1) at (0.8,1.4);
\node (bo1) at (1.3,0){$b^*$};

\begin{knot}[clip width=4pt]
\strand [blue,thick] (bi) to[out=north,in=south] (w) to [out=north, in=north] (e) to (bo1);
\strand[thick] (ci) to [out=north,in=south] (co);
\strand[thick] (ci1) to [out=north,in=south] (co1);
\flipcrossings{1,2}
\end{knot}
\end{tikzpicture}
}.$$
Here we used that $\coev_b=\ev_b^{-1}$, as $\cat{B}$ is spherical and pointed.

This $\tau$ is natural in $b$ as $\beta$ is. To see that $\tau$ is monoidal, observe that Equation \eqref{lifthalfbraidcondb} gives
\begin{eqnarray*}
\tau_{bb',c}=&\hbox{
\begin{tikzpicture}[baseline=(current  bounding  box.center)]
\node (bip) at (0,0){$b'$};
\node (bi) at (-0.3,0){$b$};
\node (ci) at (0.5,0){$c$};
\coordinate (co) at (0.5,1.8);
\coordinate (wp) at (0,0.3);
\coordinate (ep) at (1,0.6);
\coordinate (w) at (-0.3,0.5);
\coordinate (e) at (1.5, 0.9);
\node (bop) at (1,0){$(b')^*$};
\node (bo) at (1.5,0){$b^*$};

\begin{knot}[clip width=4pt]
\strand [blue,thick] (bi) to[out=north,in=south] (w) to [out=north, in=north] (e) to (bo);
\strand [blue,thick] (bip) to[out=north,in=south] (wp) to [out=north, in=north] (ep) to (bop);
\strand[thick] (ci) to [out=north,in=south] (co);
\flipcrossings{1,2}
\end{knot}
\end{tikzpicture}
}
\\
=&\tau_{b,c}\circ \tau_{b',c}.
\end{eqnarray*}
This shows that a central lift gives rise to a trivialisation of the adjoint action.

For the converse we define $\beta_b$ for $b$ a simple object by
\begin{equation}\label{halfbraidfromadtriv}
\beta_{b,c}= \tau_{b,c}\circ (\id_{bc} \otimes \coev_b),
\end{equation}
that is
$$
\hbox{
\begin{tikzpicture}[baseline=(current  bounding  box.center)]
\node (bi) at (0,0){$b$};
\node (ci) at (0.5,0){$c$};
\coordinate (co) at (0,1);
\coordinate (bo) at (0.5,1);

\begin{knot}[clip width=4pt]
\strand [blue,thick] (bi) to[out=north,in=south] (bo);
\strand[thick] (ci) to [out=north,in=south] (co);
\flipcrossings{1}
\end{knot}
\end{tikzpicture}
}=
\hbox{
\begin{tikzpicture}[baseline=(current  bounding  box.center)]
\node (bi) at (0,0){$b$};
\node (ci) at (0.5,0){$c$};
\coordinate (co) at (0.5,1.4);
\coordinate (w) at (0,0.6);
\coordinate (e) at (1,0.6);
\coordinate (c) at (1.1,0.3);
\coordinate (bo) at (1.2,1.4);

\begin{knot}[clip width=4pt]
\strand [blue,thick] (bi) to[out=north,in=south] (w) to [out=north, in=north] (e) to [out=south,in=west] (c) to [out=east,in=south] (bo);
\strand[thick] (ci) to [out=north,in=south] (co);
\flipcrossings{1}
\end{knot}
\end{tikzpicture}
}.
$$
This determines $\beta_{b,c}$ for all $b\in \cat{B}$. Similar arguments to the ones above show that this indeed satisfies the conditions for an unbraided central lift. 
\end{proof}

\subsubsection{$\cat{V}_d\hookrightarrow \cat{P}_d$ and a central lift for it}
In this subsection we will show that $\cat{V}_d\hookrightarrow \cat{P}_d$ admits a central lift. We will do this through Lemma \ref{centralliftvstriv}, using the results from Section \ref{fusionsection} to produce the trivialisation for the adjoint action $\Ad$ of $\cat{V}_d$ on $\cat{P}_d$. That is:
\begin{prop}\label{trivialisationonpd}
	The natural isomorphism $\tau\colon \Ad \Rightarrow \Id_{\cat{P}_d}$ with components for $a\in \Zd$ and $P_{m;l}\in \cat{P}_d$ given by
	$$
	\tau_{a,(m;l)}= \eta^{ am} \hbox{
		\begin{tikzpicture}[baseline=(current  bounding  box.center)]
			\node (a) at (-0.75,0){$\underline{a}$};
			\node (b) at (0.75,0){$\underline{\mi a}$};
			\node (p) at (0,0){$P_{m;l}$};
			\coordinate (v) at (0,0.6);
			\coordinate (v1) at (0,0.8);
			\node (pm) at (0,1.3){$P_{m;l}$};
			
			\draw[thick,blue] (a) to [out=north,in=west] (v);
			\draw[thick,blue] (b) to [out=north,in=east] (v1);
			\draw[thick] (p) to (v) to (v1) to (pm);
			\draw[fill,blue] (v) circle [radius=0.06];
			\draw[fill,blue] (v1) circle [radius=0.06];
		\end{tikzpicture}
	}
	$$
	is monoidal in both $a$ and $P_{m;l}$. 
\end{prop}

We remark that we want $\tau$ to be a natural isomorphism in $\Fun_\otimes (\Zd ,\Aut_\otimes(P_d))$. As $\Zd$ is a discrete category and $\cat{P}_d$ is semi-simple, showing naturality is trivial. We need to show that it is monoidal on $\Zd$, and monoidal on $\cat{P}_d$. Recall that $\Ad_{a}\circ\Ad_b \cong \Ad_{a+b}$ by fusion of $\underline{a}$ and $\underline{b}$ along $\psi^L_{a;(b;0)}$ and simultaneously fusing $\underline{\mi a}$ and $\underline{\mi b}$ to $\underline{\mi a \mi b}$, so the monoidality on $\Zd$ is:

\begin{lem}
Let $a,b\in \Zd$ and $P_{m;l}\in \cat{P}_d$. Then
$$
\eta^{(a+b)m}
\hbox{
\begin{tikzpicture}[baseline=(current  bounding  box.center)]

\node (a) at (-1,0){$\underline{a}$};
\node (b) at (-0.6,0){$\underline{b}$};
\node (p) at (0,0){$P_{m;l}$};
\coordinate (v) at (0,0.7);
\coordinate (v1) at (0,1);
\coordinate (pm) at (0,1.6);

\node (am) at (1,0){$\underline{\mi a}$};
\node (bm) at (0.6,0){$\underline{\mi b}$};
\coordinate (vm) at (0,0.8);
\coordinate (v1m) at (0,1.1);

\draw[thick,blue] (am) to [out=north,in=east] (v1m);
\draw[thick,blue] (bm) to [out=north,in=east] (vm);
\draw[fill,blue] (vm) circle [radius=0.06];
\draw[fill,blue] (v1m) circle [radius=0.06];

\draw[thick,blue] (a) to [out=north,in=west] (v1);
\draw[thick,blue] (b) to [out=north,in=west] (v);
\draw[thick] (p) to (v) to (v1) to (pm);
\draw[fill,blue] (v) circle [radius=0.06];
\draw[fill,blue] (v1) circle [radius=0.06];
\end{tikzpicture}
}
=
\eta^{(a+b)m}
\hbox{
\begin{tikzpicture}[baseline=(current  bounding  box.center)]

\node (a) at (-1,0){$\underline{a}$};
\node (b) at (-0.6,0){$\underline{b}$};
\node (p) at (0,0){$P_{m;l}$};
\coordinate (v) at (-0.6,0.6);
\coordinate (v1) at (0,1);
\node (pm) at (0,1.6){};

\node (am) at (1,0){$\underline{\mi a}$};
\node (bm) at (0.6,0){$\underline{\mi b}$};
\coordinate (vm) at (1,0.6);
\coordinate (v1m) at (0,1.1);

\draw[thick,blue] (am) to [out=north,in=south] (vm)to [out=north, in=east] (v1m);
\draw[thick,blue] (bm) to [out=north,in=west] (vm) ;
\draw[fill,blue] (vm) circle [radius=0.06];
\draw[fill,blue] (v1m) circle [radius=0.06];

\draw[thick,blue] (a) to [out=north,in=west] (v);
\draw[thick,blue] (b) to [out=north,in=south] (v) to [out=north, in=west] (v1);
\draw[thick] (p) to (v1) to (pm);
\draw[fill,blue] (v) circle [radius=0.06];
\draw[fill,blue] (v1) circle [radius=0.06];
\end{tikzpicture}
}.
$$
\end{lem}
\begin{proof}
Apply Lemma \ref{psileftrightorder}, then Lemmas \ref{leftpsiasso} and \ref{rightpsiasso}, followed by Lemma \ref{psisymmetry}.
\end{proof}

Recall that the isomorphism $\Ad_{a}(P P') \cong  \Ad_{a}(P) \otimes \Ad_a(P')$ is given by $\id_{\underline{a}P}\otimes\coev_{-a}\otimes\id_{\underline{\mi a}P'}$. The monoidality in $\mathcal{P}_{d}$ is then a consequence of Proposition \ref{conjugationtrivfusion}:
\begin{cor}
	Let $P_{m;l},P_{m';l'},P_{m'';l''}\in \mathcal{P}_d$ and $a\in \Zd$. Then
	$$
	\eta^{am''}
	\hbox{
		\begin{tikzpicture}[baseline=(current  bounding  box.center)]
			\node (a) at (-1.4,-2){$\underline{a}$};
			\node (am) at (1.4,-2){$\underline{\mi a}$};
			\node (pl) at (-0.5,-2){$P_{m;l}$};
			\node (pr) at (0.5,-2){$P_{m';l'}$};
			\coordinate (v1) at (0,-1.2);
			\coordinate (v2) at (0,-0.8);
			\coordinate (v3) at (0,-0.6);
			\node (pm) at (0,0){$P_{m'';l''}$};
			
			\draw[thick,blue] (a) to [out=north,in=west] (v2);
			\draw[thick,blue] (am) to [out=north,in=east] (v3);
			\draw[thick] (pl) to  [out=north, in=west] (v1) to [out=east, in=north] (pr);
			\draw[thick] (pm) to (v1);
			\draw[fill] (v1) circle [radius=0.05];
			\draw[fill,blue] (v2) circle [radius=0.05];
			\draw[fill,blue] (v3) circle [radius=0.05];
		\end{tikzpicture}
	}
	=
	\eta^{a(m+m')}
	\hbox{
		\begin{tikzpicture}[baseline=(current  bounding  box.center)]
			\node (a) at (-1.4,-2){$\underline{a}$};
			\node (am) at (1.4,-2){$\underline{\mi a}$};
			\node (pl) at (-0.5,-2){$P_{m;l}$};
			\node (pr) at (0.5,-2){$P_{m';l'}$};
			\coordinate (v1) at (0,-0.7);
			\coordinate (v2l) at (-0.5,-1);
			\coordinate (v2r) at (0.5,-1);
			\coordinate (v3l) at (-0.5,-1.2);
			\coordinate (v3r) at (0.5,-1.2);
			\coordinate (c) at (0,-1.4);
			\node (pm) at (0,0){$P_{m'';l''}$};
			
			\draw[thick,blue] (a) to [out=north,in=west] (v3l);
			\draw[thick,blue] (am) to [out=north,in=east] (v3r);
			\draw[thick] (pl) to (v3l) to [out= north, in=south] (v2l) to [out=north, in=west] (v1) to [out=east, in=north] (v2r) to  (v3r) to (pr);
			\draw[thick] (pm) to (v1);
			\draw[thick,blue] (v2l) to [out =east, in=west] (c) to [out=east, in=west] (v2r);
			\draw[fill] (v1) circle [radius=0.06];
			\draw[fill,blue] (v2l) circle [radius=0.06];
			\draw[fill,blue] (v2r) circle [radius=0.06];
			\draw[fill,blue] (v3r) circle [radius=0.06];
			\draw[fill,blue] (v3l) circle [radius=0.06];
		\end{tikzpicture}
	}.
	$$
\end{cor}
\begin{proof}
	From the fusion rules Equation \eqref{fusionrulesPd} we see that $m''=m+m'+\frac{1}{2}(l+l'-\nu)$, with $\nu$ the summation variable, this gives the exponent of $\eta^a$ on the left hand side. Using Proposition \ref{conjugationtrivfusion}, we see that the exponent of $\eta^a$ on the right hand side is equal to this.
\end{proof}

This concludes the proof of Proposition \ref{trivialisationonpd}. Using Lemma \ref{centralliftvstriv} now gives an unbraided central lift for $\cat{V}_d \subset \mathcal{P}_d$. We will now show that this unbraided lift is in fact braided, for the braiding on $\cat{V}_d$ as induced from the quadratic form $q_d$ (Definition \ref{vddef}). It suffices to show that the self-braiding on $\underline{a}$ from the central lift (Equation \eqref{halfbraidfromadtriv}) agrees with the self-braiding on $\cat{V}_d$, which is given by $q_{d}(a) \id_{\underline{a}\otimes\underline{a}}=\eta^{a^2} \id_{\underline{a}\otimes\underline{a}}$. 

\begin{lem}
	For $a\in \Zd$ we have
	$$
	\eta^{a^2}
	\hbox{
		\begin{tikzpicture}[baseline=(current  bounding  box.center)]
			\node (a) at (-0.75,0){$\underline{a}$};
			\coordinate (b) at (0.7,0);
			\node (aout) at (0.9,1.3){$\underline{a}$};
			\node (p) at (0,0){$\underline{a}$};
			\coordinate (v) at (0,0.6);
			\coordinate (v1) at (0,0.8);
			\node (pm) at (0,1.3){$\underline{a}$};
			
			\draw[thick,blue] (a) to [out=north,in=west] (v);
			\draw[thick,blue] (aout) to [out=south, in=east] (b) to [out=west,in=east] (v1);
			\draw[thick,blue] (p) to (v) to (v1) to (pm);
			\draw[fill,blue] (v) circle [radius=0.06];
			\draw[fill,blue] (v1) circle [radius=0.06];
		\end{tikzpicture}
	}
	=\eta^{a^2}
	\hbox{
		\begin{tikzpicture}[baseline=(current  bounding  box.center)]
			\node (a) at (-0.75,0){$\underline{a}$};
			\node (aout) at (-0.75,1.3){$\underline{a}$};
			\node (p) at (0,0){$\underline{a}$};
			\node (pm) at (0,1.3){$\underline{a}$};
			
			\draw[thick,blue] (a) to [out=north,in=south] (aout);
			\draw[thick,blue] (p) to (pm);
		\end{tikzpicture}
	}.
	$$
\end{lem}
\begin{proof}
	By Lemma \ref{psileftrightorder}, we can swap the order of the fusions on the left hand side. The first fusion is then 
	$$
	\psi^R_{(a;0),\mi a}=\psi^L_{a, (\mi a,0)}=\ev_{\mi a},
	$$
	by Lemmas \ref{psisymmetry} and \ref{psievagree}. Now applying the snake relation for $\ev_{\mi a}$ and $\coev_{-a}$ gives the result.
\end{proof}

In summary:

\begin{thm}\label{pdisvdmod}
	The category $\mathcal{P}_d$ can be given the structure of a free spherical $\mathcal{V}_d$-module fusion category. 
\end{thm}

\section{The tensor equivalence}\label{Section4}

We are now in a position to prove our main result. Let us summarise what we have found so far. We have established that both $\cns$ (Proposition \ref{cnsisvdmod}) and $\cat{P}_d$ (Theorem \ref{pdisvdmod}) are spherical free $\cat{V}_d$-module fusion categories. We have further established that $\cns$ is generated by the charge $-1$ self-dual object $[1,1]$ of quantum dimension $2 \cos (\frac{\pi}{d})$, Theorem \ref{cnsgenerated}. Recall that in the spherical pivotal structure $P_{0;1}$ also has quantum dimension $2 \cos (\frac{\pi}{d})$ (Equation \eqref{quantumdimensionofP01}. Furthermore, its dual is $P_{-1,1}$, making it self-dual of charge $-1$. This gives a tensor functor $\cns\rar \pd$ mapping $[1,1]$ to $P_{0;1}$. Using semi-simplicity one sees that this functor is an equivalence. That is, we have:

\begin{thm}\label{badassmaintheorem}
The categories $\cns$ and $\cat{P}_d$ are equivalent $\vd$-module tensor categories for any $d$, with equivalence that sends $[l,2m+l] \mapsto P_{m;l}$.
\end{thm}

In particular, $\cns$ and $\pd$ are equivalent as spherical fusion categories.

\begin{rmk}
Following on from Remark \ref{drcrcompareremark1}, it is educative to compare our methods with \cite{Davydov2018}. That remark addresses how our view on $\cns$ changes. What it does not address is the absence of equivariant duality data for $P_{0;1}$. To recover this, note that we have provided something stronger: we have a trivialisation of the adjoint action of $\cat{V}_d$ on all of $\cat{P}_d$. Specialising this trivialisation to $P_{0;1}$ gives rise to the data used in \cite{Davydov2018}.
\end{rmk}

\section*{Acknowledgements}
The authors gratefully acknowledge the Topology and Geometry conference (2015) and the Blow Up in Bonn for making sure our paths crossed. ARC is supported by Cardiff University. We also gratefully acknowledge support from Cardiff University for hospitality for Thomas during his visit. TAW is supported by Lincoln College, University of Oxford. TAW is grateful to Andr\'e Henriques, M\'arton Habliscek, and Christoph Weiss for helpful discussions during the preparation of this paper. TAW also wants to thank the wardens of Skokholm Island, the Oxford Zoology Knowles Lab, and autumnal storms for a long and beautiful stay on an isolated island that inspired this paper.

\newcommand{\etalchar}[1]{$^{#1}$}

\end{document}